\documentclass[11pt,a4paper, notitlepage]{article}
\usepackage{import} 
\usepackage[tmargin=2.5cm,bmargin=4cm,lmargin=1.25cm,rmargin=1.25cm]{geometry}
\usepackage[utf8]{inputenc}
\usepackage[english]{babel}

\usepackage[T1]{fontenc} 
\usepackage{amssymb,amsfonts}
\usepackage{amsthm}
\usepackage{tikz-cd}
\usepackage[colorinlistoftodos]{todonotes}
\usepackage[framemethod=TikZ]{mdframed} 
\usepackage{calligra,mathrsfs}
\usepackage{amsmath}
\usepackage{mathtools}
\usetikzlibrary{shapes.geometric, arrows}
\usepackage[backend=biber,style=alphabetic,sorting=nyt, maxnames=10, maxalphanames=10]{biblatex}
\usepackage{csquotes}
\usepackage{colonequals}
\usepackage{adjustbox} 
\usepackage[colorlinks, final]{hyperref} 
\usepackage{cleveref}  
\usepackage{enumerate} 
\definecolor{imperialred}{RGB}{237, 41, 57}
\definecolor{royalblue}{RGB}{64, 106, 212}
\definecolor{link}{RGB}{11,0,128}
\hypersetup{
	colorlinks = true,
	breaklinks = true,
	citecolor=royalblue,
	linkcolor=imperialred,
	urlcolor=link,
	linktocpage = true
}


\newrobustcmd*{\citeauthorfullname}{\AtNextCite{\DeclareNameAlias{labelname}{given-family}}\citeauthor}

\newtheorem{thm}{Theorem}[section]
\newtheorem{cor}[thm]{Corollary}
\newtheorem{prop}[thm]{Proposition}
\newtheorem{lem}[thm]{Lemma}
\newtheorem{conj}[thm]{Conjecture}

\newtheorem*{quest*}{Question}

\theoremstyle{definition}
\newtheorem{defn}[thm]{Definition}

\theoremstyle{remark}

\newtheorem{rem}[thm]{Remark}
\newtheorem{rems}[thm]{Remarks}

\newtheoremstyle{montobbo}
{7pt}
{5pt}
{}
{}%
{\bfseries}
{}%
{.5em}
{\thmnumber{\@{#1}{}~\@{#2}.}%
	\thmnote{~{\bfseries#3.}}}

\theoremstyle{montobbo}
	\newtheorem{montobo}[thm]{}
	\newcommand{\bmo}{\begin{montobo}}
	\newcommand{\emo}{\end{montobo}}
	\newtheorem{talk}{Lecture}
	\newcommand{\btalk}{\begin{talk}}
		\newcommand{\etalk}{\end{talk}}
	
\newenvironment{manualtheorem}[1]{
	\renewcommand\thethm{#1}
	\thm
}{\endthm}
\newcommand{\but}[1]{\begin{manualtheorem}{#1}}
\newcommand{\eut}{\end{manualtheorem}}
\newenvironment{manualproposition}[1]{
	\renewcommand\thethm{#1}
	\prop
}{\endprop}
\newcommand{\bup}[1]{\begin{manualproposition}{#1}}
	\newcommand{\eup}{\end{manualproposition}}
\numberwithin{equation}{thm}
\newenvironment{teqn}{\begin{equation}\begin{tikzcd}\displaystyle}{\end{tikzcd}\end{equation}}
\newenvironment{teqn*}{\begin{equation*}\begin{tikzcd}\displaystyle}{\end{tikzcd}\end{equation*}}
\newenvironment{tpic}{\begin{equation}\begin{tikzpicture}}{\end{tikzpicture}\end{equation}}
\newenvironment{tpic*}{\begin{equation*}\begin{tikzpicture}}{\end{tikzpicture}\end{equation*}}

\DeclareMathAlphabet\matheu{U}{eus}{m}{n}
\DeclareMathAlphabet{\mathmm}{U}{mmambb}{m}{n}
\newcommand{\bb}[1]{\mathbb{#1}}
\newcommand{\fr}[1]{\mathfrak{#1}}
\newcommand{\ca}[1]{\mathcal{#1}}
\newcommand{\scr}[1]{\mathscr{#1}}

\newcommand{\comment}[1]{}
\newcommand{\bun}{\begin{itemize}}
\newcommand{\bn}{\begin{enumerate}}
\newcommand{\bt}{\begin{thm}}
\newcommand{\bl}{\begin{lem}}
\newcommand{\bp}{\begin{prop}}
\newcommand{\bc}{\begin{cor}}
\newcommand{\bd}{\begin{teqn}}
\newcommand{\bud}{\begin{teqn*}\textstyle}	
\newcommand{\bs}{\begin{proof}}
\newcommand{\br}{\begin{rem}}
\newcommand{\bdf}{\begin{defn}}
\newcommand{\bcj}{\begin{conj}}

\newcommand{\eun}{\end{itemize}}
\newcommand{\en}{\end{enumerate}}
\newcommand{\et}{\end{thm}}
\newcommand{\el}{\end{lem}}
\newcommand{\ep}{\end{prop}}
\newcommand{\ec}{\end{cor}}
\newcommand{\ed}{\end{teqn}}
\newcommand{\eud}{\end{teqn*}}
\newcommand{\es}{\end{proof}}
\newcommand{\er}{\end{rem}}
\newcommand{\edf}{\end{defn}}
\newcommand{\ecj}{\end{conj}}
\DeclareMathOperator{\etale}{\acute{e}t}

\DeclareMathOperator{\id}{id}
\DeclareMathOperator{\Hom}{Hom}

\DeclareMathOperator{\dep}{depth}
\DeclareMathOperator{\trdeg}{tr.deg}

\DeclareMathOperator{\Ext}{Ext}

\DeclareMathOperator{\depth}{depth}
\DeclareMathOperator{\codim}{codim}

\DeclareMathOperator{\Gal}{Gal}

\DeclareMathOperator{\colim}{colim}
\newcommand{\units}[1]{#1^{\times}}

\DeclareMathOperator{\Frac}{Frac}
\DeclareMathOperator{\spec}{Spec}
\DeclareMathOperator{\mspec}{MaxSpec}

\DeclareMathOperator{\Pic}{Pic}

\newcommand{\tir}{\arrow[r]}

\newcommand{\iso}{\xrightarrow{\sim}}
\newcommand{\isom}{\arrow[r, "\sim"]}



\bibliography{/Users/arnabkundu/Preamble/bibli.bib} 
\makeatletter                                                                                                                             

\newtheorem{pl}[thm]{Presentation Lemma }
\newenvironment{claim}[1]{\par\comment{\noindent}\underline{Claim}#1:}{}
\newenvironment{whatever}[2]{\par\comment{\noindent}\underline{#1}#2:}{}
\newenvironment{claimproof}[1]{\par\comment{\noindent}\underline{Proof:}\space#1}{}

\newcommand{\ccite}[2]{\cite[#2]{#1}}
\newcommand{\stacks}[1]{\ccite{stacks-project}{\href{https://stacks.math.columbia.edu/tag/#1}{Tag \MakeUppercase{#1}}}}

\newcommand{\inv}[1]{[\frac{1}{#1}]}
\DeclareMathOperator{\coh}{Coh}
\DeclareMathOperator{\rings}{CRings}
	
\title{\uppercase{\textbf{\large{Gersten's Injectivity for Smooth Algebras over Valuation Rings}}}}
\author{\uppercase{Arnab Kundu}{\let\thefootnote\relax\footnote{University of Toronto, St.~George Campus, Toronto, Canada\newline\hspace*{0.48cm} Email: arnab.kundu@utoronto.ca.}}}
\date{}
\usepackage{doi}

\begin{document}
\maketitle
\abstract{Gersten's injectivity conjecture for a functor $F$ of ``motivic type'', predicts that given a semilocal, ``non-singular'', integral domain $R$ with a fraction field $K$, the restriction morphism induces an injection of $F(R)$ inside $F(K)$. We prove two new cases of this conjecture for smooth algebras over valuation rings. Namely, we show that the higher algebraic $K$-groups of a semilocal, integral domain that is an essentially smooth algebra over an equicharacteristic valuation ring inject inside the same of its fraction field. Secondly, we show that Gersten's injectivity is true for smooth algebras over, possibly of mixed-characteristic, valuation rings in the case of torsors under tori and also in the case of the Brauer group.}
\tableofcontents
{\let\thefootnote\relax\footnote{April 2024}}
\section{Gersten's Injectivity Conjecture}
	 Let $\cal{C}$ be an additive category, for example, the category of Abelian groups, and let $\cal{S}\subseteq\rings$ be a full subcategory of the category of commutative, unital rings. Given a functor $\scr{F}\colon\rings\to\cal{C}$, we say that $\scr{F}$ \textit{satisfies Gersten's injectivity for $\cal{S}$}, if for any semilocal, integral domain $R\in\cal{S}$ with a fraction field $K$, we have an injection
 	\bd\label{diag:injectivity}\scr{F}(R)\arrow[r, hook]&\scr{F}(K).\ed Gersten's injectivity conjecture imprecisely stated predicts the following.
	\begin{conj}\label{conj:Gersten_injectivity} A functor $\scr{F}$ of `motivic type' satisfies Gersten's injectivity for a subcategory $\cal{S}$ of `non-singular' rings. 
	\end{conj}
	Traditionally, for the subcategory $\mathrm{Reg}$ of regular rings, we expect \Cref{conj:Gersten_injectivity} to be true for the functor of algebraic $K$-groups, the Milnor $K$-groups, the Hermitian Witt groups, the de Rham cohomology groups, the \'etale cohomology groups with coefficients in the $\ell$-th roots of unity, the motivic cohomology groups, etc. Validity of a certain `effacement theorem' (in the style of \ccite{ct_hoobler_kahn}{Theorem~2.2.7}) is at the foundation of \Cref{conj:Gersten_injectivity} for each of above mentioned functors. Knowing that a functor $\scr{F}$ satisfies Gersten's injectivity has diverse practical benefits, in particular, the behaviour of such a functor $\scr{F}$ is influenced by its values on fields. To provide background for our discussion, let's briefly delve into Gersten's conjecture.
	\subsubsection*{Digression: Gersten's conjecture}
		Let $H^i(-)$ denote one of the cohomology theories mentioned above. Gersten's injectivity forms a valuable part of Gersten's conjecture, which, in essence, predicts that for a regular, Noetherian scheme $X$, the group $H^n(X)$ can be read off by calculating $H^{n-i}(\kappa(x))$, where $0\le i\le n$ and $\kappa(x)$ is the residue field of $X$ at a point $x\in X$ of codimension $i$. Gersten's conjecture has been extensively studied in the literature thanks to its far reaching consequences and wide range of applications. Let us note some of them below.
		\bun
			\item[$\circ$] In the case of algebraic $K$-theory, Gersten's conjecture illuminates the connection of the algebraic cycles with algebraic $K$-theory (see Bloch's formula proved in \ccite{quillen_k_theory}{\textsection7, Theorem~5.19}). Using Bloch's formula and the cup product structure on algebraic $K$-groups, one can readily define intersection products of algebraic cycles in regular, Noetherian schemes (see \cite{grayson_products_K_theory_and_intersecting_algebraic_cycles}).
			\item[$\circ$] In the case of Milnor $K$-theory, Gersten's conjecture paves the way to prove Levine's generalised Bloch--Kato conjecture for semi-local, equicharacteristic rings. Furthermore, one can deduce Beilinson's conjecture from Gersten's conjecture (see \cite{kerz_gersten_milnor_k-theory}, which simultaneously proves Gersten's conjecture in the equicharacteristic case).
			\item[$\circ$] In the case of de Rham cohomology, Gersten's conjecture can be used to verify Washnitzer's conjecture relating the filtration by coniveau with the same by hypercohomology (see \ccite{bloch_ogus_gersten}{Corollary~6.9}).
		\eun
		Therefore, unsurprisingly, Gersten's conjecture enjoys a pivotal role in the study of each of these functors. 
	 \vspace{0.25cm}\newline Warmed up by the above digression, we turn our attention back to \Cref{conj:Gersten_injectivity}. Given the significance of the conjecture, it seems completely justified to question its purview. Hence, we may pose the following. 
	 \begin{quest*}
	 	Is there a larger category $\cal{S}\supset\mathrm{Reg}$ of `non-singular' rings for which we expect Gersten's injectivity to be satisfied by a functor $\scr{F}$ of `motivic type'?
	 \end{quest*} 
	 As a matter of fact, the answer is yes. Indeed, thanks to Zariski's local uniformisation conjecture (see, for example, \ccite{phd-thesis}{Conjecture~2.1.1}), which predicts that any valuation ring is ind-regular, and a limit argument, any functor that commutes with filtered colimits and satisfies \Cref{conj:Gersten_injectivity} for $\mathrm{Reg}$, also does the same for the category $\cal{S}_{\mathrm{Val}}$ of essentially smooth algebras over valuation rings. However, as mentioned in loc.~cit., although expected to be true, Zariski's conjecture is widely open. This leaves the possibility for a challenge to prove \Cref{conj:Gersten_injectivity} unconditionally for $\cal{S}_{\mathrm{Val}}$ in the case of functors $\scr{F}$ for which \Cref{conj:Gersten_injectivity} is known to be true for a subcategory of $\mathrm{Reg}$. 
	\par The purpose of this article is twofold. We prove \Cref{conj:Gersten_injectivity} in the case of two different classes of functors, namely
	\bun
		\item[$\circ$] the higher algebraic $K$-groups, and 
		\item[$\circ$] the functor that classifies $T$-torsors and the functor that classifies $T$-gerbes, where $T$ is some fixed torus.
	\eun
	We write the precise statements of the two theorems below before having separate discussions on each of them later. Let $V$ be a valuation ring and let $R$ be a semilocalisation of a smooth $V$-algebra with a fraction field $F$.
	\bt[\Cref{cor:Gersten_injectivity_for_valuation}]\label{thm:1.1}
		If $V$ contains a field, then we have an injection \bud K_i(R)\arrow[r,hook]& K_i(F), \hspace{0.25cm}\text{for all }i.\eud
	\et
	\bt[\Cref{cor:torus_Gersten_injectivity}]\label{thm:1.2}
		Given an $R$-torus $T$, we have injections \bud H^1_{\etale}(R, T)\arrow[r,hook]& H^1_{\etale}(F, T) \hspace{0.5cm}\text{ and }\hspace{0.5cm}H^2_{\etale}(R, T)\arrow[r,hook]& H^2_{\etale}(F, T).\eud
	\et
	Plugging $T=\bb{G}_m$ in the displayed formula on the right, we get Gersten's injectivity for $\cal{S}_{\mathrm{Val}}$ in the case of Brauer groups.
	\subsection*{Background and Remarks on \Cref{thm:1.1}}
	Let us start this subsection with a review of the known cases of Gersten's conjecture (see \ccite{quillen_k_theory}{\textsection7, Conjecture~5.10}) in the case of algebraic $K$-theory.
	\bun
		\item[$\circ$] \citeauthor{gersten} in \ccite{gersten}{Theorem~1.3} used representation theoretic techniques to prove Gersten's conjecture for discrete valuation rings with a finite residue field. \citeauthor{sherman_group_representation} in \cite{sherman_group_representation} generalised Gersten's method to prove the same for discrete valuation rings whose residue field is an algebraic extension of a finite field.
		\item[$\circ$] \citeauthor{quillen_k_theory} in \ccite{quillen_k_theory}{\textsection 7, Theorem~5.11} proved Gersten's conjecture for essentially smooth algebras over a field. The key input is a certain `geometric presentation lemma', which will appear in our discussion later.
		\item[$\circ$] \citeauthor{gillet_levine} in \ccite{gillet_levine}{Corollary~6} subsequently generalised Quillen's method to the mixed-characteristic. They showed that Gersten's conjecture is true for all essentially smooth algebras over discrete valuation rings if the same holds for all discrete valuation rings.
		\item[$\circ$] \citeauthor{kelly-morrow} in \ccite{kelly-morrow}{Theorem~3.1} recently proved Gersten's injectivity for equicharacteristic valuation rings. Their method is to pass to special classes of valuation rings that we already know are ind-regular, so as to reduce to \ccite{quillen_k_theory}{\textsection 7, Theorem~5.11}.
	\eun
	\par In \textsection\ref{section:G-theory}, our goal is to establish \Cref{thm:1.1}. The technique of the proof is based on the same of \ccite{gillet_levine}{Theorem}. More precisely, we prove the following generalisation of loc.~cit. Let $V$ be a valuation ring and let $R$ be a semilocalisation of a smooth $V$-algebra. Suppose that $\ca{P}\subset\spec(R)$ is the subset that corresponds to the generic points of the $V$-special fibre of $\spec(R)$. Let us denote the semilocalisation of $R$ at the primes in $\ca{P}$ by $R_{\ca{P}}$.
	\bt[\Cref{thm:4}]\label{thm:1.3}
		 The canonical restriction induces an injection \bud K_i(R)\arrow[r,hook]& K_i(R_{\ca{P}}), \hspace{0.25cm}\text{for all }i.\eud
	\et
	Thanks to \Cref{lem:snail_generic_point_of_special_fibre}, the semilocal ring $R_{\ca{P}}$ is a Pr\"ufer domain (i.e., it's an integral domain whose local rings are valuation rings). As a consequence, \Cref{thm:1.3} reduces the verification of Gersten's injectivity for $\ca{S}_{\mathrm{val}}$ to that of semilocal Pr\"ufer domains. Finally, we deduce \Cref{thm:1.1} from \ccite{kelly-morrow}{Theorem~3.1}.
	\par It is worth noting that, unlike in \Cref{thm:1.1}, the ring $V$ in \Cref{thm:1.3} is not assumed to be equicharacteristic. The main technical component of the proof of this theorem is Presentation Lemma~\ref{lem:presentation_lemma_Prufer}, which is a ``geometric presentation lemma'' in the style of Quillen in \ccite{quillen_k_theory}{\textsection7, Lemma~5.12} (see \ccite{ces_torsors}{beginning of \textsection4.1}). \citeauthor{gillet_levine} in \ccite{gillet_levine}{Lemma~1} proved a version of the presentation lemma over mixed-characteristic, discrete valuation rings (cf.~\ccite{luders_relative_gersten_conjecture_for_milnor_k_theory}{Lemma~2.12} and \ccite{ces_grothendieck-serre}{Variant~3.7}). Employing techniques from the proofs of \ccite{ces_grothendieck-serre}{Variant~3.7} and \ccite{phd-thesis}{Proposition~6.4}, we construct our version of the presentation lemma (\Cref{lem:presentation_lemma_Prufer}) over valuation rings of finite Krull dimension, generalising \ccite{gillet_levine}{Lemma~1}.
	\subsection*{Background and Remarks on \Cref{thm:1.2}}
	Before delving into the discussion on \Cref{thm:1.2}, let us contextualise it. \Cref{conj:Gersten_injectivity} in the case of the functors $H^1(-,\bb{G}_m)=\Pic(-)$ and $H^2(-,\bb{G}_m)$ for regular rings was explored by Grothendieck in his seminal paper \cite{brauerII}. He presented an \'etale cohomological interpretation of the well-known Brauer groups, ultimately establishing Gersten's injectivity by analysing the long exact sequence associated to a certain `divisor short exact sequence'. Subsequently, \citeauthor{ct-sansuc_flasque_tori}~\cite{ct-sansuc_flasque_tori} generalised these results to the functors defined by non-split tori, establishing Gersten's injectivity for them. Their key contribution lies in the introduction of the concept of a particular class of isotrivial tori termed `flasque tori', characterised by simple Galois theoretic data. Remarkably, they demonstrated that any isotrivial tori can be resolved by flasque tori. Recently, employing similar techniques, Guo~\cite{ning_valuation} proved Gersten's injectivity for valuation rings.
	\par Inspired by the established techniques, in \textsection\ref{section:toral_case}, our objective is to prove \Cref{thm:1.2}. The strategy is roughly based on the proof of \ccite{ct-sansuc_flasque_tori}{Theorem~2.2}. Considering an $R$-torus $T$, the key steps are as follows: 
	\bun
		\item[$\circ$] Firstly, leveraging the fact that $T$ has a flasque resolution, denoted by $E^{-}$, (refer to \textsection\ref{montobo:flasque_torus_resolution}), an analysis of the long exact sequence of cohomology associated to $E^{-}$ reduces to show \Cref{thm:1.2} in the case when $T$ is flasque. 
		\item[$\circ$] Thereafter, writing down the long exact sequence of cohomology with supports in a closed subscheme and subsequently, applying the coniveau spectral sequence associated to the filtration by supports, we further reduce to show the following vanishing statement of local cohomology in low degrees.
		\bt\label{thm:1.5}
			Let $\frak{p}\subset R$ be a prime ideal. For a flasque $R$-torus $T$, there is vanishing \bud H^{q}_{\frak{p}}(R, T)=0, \hspace{1cm}\text{for }q\le 2.\eud
		\et
		\item[$\circ$] Finally, we establish purity for torsors under tori (see \Cref{prop:extend_tori_torsor}), which we derive as a consequence of a weak version of the Auslander--Buchsbaum formula for smooth algebras over valuation rings (see \Cref{lem:Gabber_Ramero_S2_condition}).
	\eun
	\par The results of this section appeared in the author's thesis \ccite{phd-thesis}{Chapter~4} and they were simultaneously and independently obtained by \citeauthor{ning_liu_quasi-split} in \cite{ning_liu_quasi-split}.
	\subsection*{Notations and conventions}
		Let $I$ be an ideal in a ring $A$.
		\bun			
			\item[$\circ$] The vanishing locus of $I$ is denoted by $V(I)\subseteq\spec A$. If $I$ is principal with generator, say $t$, then $V(t)$ denotes $V(I)$.
			\item[$\circ$] The maximal spectrum of $A$ is denoted by $\mathrm{MaxSpec}(A)\subseteq\spec A$.
			\item[$\circ$] If $I$ is prime, then the localisation of $A$ at $I$ shall be denoted by $A_I$.
			\item[$\circ$] If $A$ is an integral domain, then the fraction field of $A$ is denoted by $\Frac(A)$.
			\item[$\circ$] An $A$-algebra is called \textit{essentially smooth} if it can be obtained as the semilocalisation at finitely many primes of a smooth $A$-algebra.
			\item[$\circ$] Given an $A$-scheme $X$ and an algebra $a\colon A\to A'$, the base change of $X$ along $a$ is denoted by $X_{A'}$.
			\item[$\circ$] If $I$ is prime, the residue field of $A$ at $I$ is denoted by $\kappa(I)$.
		\eun
	\subsection*{Acknowledgements}
	The author extends gratitude to Matthew Morrow and Elden Elmanto for their encouragement, which inspired the author to pursue Gersten's injectivity conjecture in the non-Noetherian case. Special thanks are also extended to Tess Bouis, Ofer Gabber, Ning Guo, Marc Levine, Morten L\"uders, David Rydh and Marco Volpe for their insightful conversations. This project started when the author was a PhD student and he is intellectually indebted to  K\k{e}stutis \v{C}esnavi\v{c}ius for his constant support throughout. This project received partial funding from the European Research Council under the European Union’s Horizon 2020 research and innovation program (grant agreement No. 851146).
\section{Generalities on Pr\"ufer Domains}\label{section:Prufer}
	In this section, we establish several foundational concepts regarding Pr\"ufer domains. These rings play a pivotal role in this article; for instance, they naturally emerge in the proof of the crucial technical foundation of \textsection\Cref{section:G-theory}, namely, Presentation Lemma~\ref{lem:presentation_lemma_Prufer}. To establish this lemma, we rely on cut-and-paste type techniques, aided by Lemmas~\ref{lem:Prufer_domain_gluing} and \ref{lem:Prufer_domain_cutting_and_pasting}, along with an approximation result, \Cref{lem:approximation:prufer_domain}. Additionally, the application of \Cref{lem:snail_generic_point_of_special_fibre} recurs throughout the article (see \Cref{thm:4} and \Cref{lem:Brauer_group_injects}).
	\par We start with the definition of Pr\"ufer domains.
	\bdf[\ccite{gilmer_multiplicative_ideal_theory}{\textsection22}]
		A \textit{Pr\"ufer domain} is an integral domain whose localisation at every prime ideal is a valuation ring.
	\edf
	For equivalent definitions of Pr\"ufer domains, see \ccite{gilmer_multiplicative_ideal_theory}{Theorem~22.1} and \stacks{092S}.
	\par We state some permanence properties of Pr\"ufer domains below.
	\bl\label{lem:pemanence_properties} 
		Given a Pr\"ufer domain $R$, a prime ideal $\frak{p}\subset R$ and a multiplicative subset $S\subset R$, the following are Pr\"ufer domains$\colon$
		\bn[$\mathrm{(}$a$\mathrm{)}$]
			\item\label{localisation:pt_b} the localisation $S^{-1}R$, and
			\item\label{quotient:pt_c} the quotient $R/\frak{p}$.
		\en
	\el
	In \Cref{lem:prufer_domain}, we will demonstrate that the concept of semilocal Pr\"ufer domains coincides with another significant class of rings known as the semilocal Krull domains (see \Cref{defn:Krull_domain} below). We will take advantage of this fact in Lemmas~\ref{lem:Prufer_domain_gluing}, \ref{lem:Prufer_domain_cutting_and_pasting} and \ref{lem:approximation:prufer_domain}. 
	\bdf[cf.~\ccite{matsumura_comm_ring_theory}{Chapter 4, \textsection12}]\label{defn:Krull_domain}
	An integral domain $R$ with a fraction field $K$ is called a \textit{semilocal Krull domain} if there exist a nonempty finite set $\Lambda$ and valuations rings $\{R_{\lambda}\}_{\lambda\in\Lambda}$ in $K$ such that $R=\bigcap_{\lambda\in\Lambda}R_{\lambda}$, where the intersection is taken in $K$.
	\edf
	\begin{rems}~
		\bn
			\item \Cref{defn:Krull_domain} extends the one in loc.~cit.~because we allow non-discrete valuations.
			\item Contrary to loc.~cit., we restrict \Cref{defn:Krull_domain} to the semilocal setting for simplicity. This particular case will be adequate for our purposes.
		\en
	\end{rems}
	\bl\label{lem:prufer_domain}
		A Pr\"ufer domain $R$ is the intersection in $\Frac(R)$ of the valuation rings obtained by the localisations at the maximal ideals. Additionally, if $R$ is semilocal, 
		\bn[$\mathrm{(}$a$\mathrm{)}$]
			\item\label{point:2_prufer_domain} for a subfield $K\subseteq\Frac(R)$, the intersection $R\cap K$ is a semilocal Pr\"ufer domain with $\Frac(R\cap K)=K$ such that the canonical morphism $\mspec(R)\twoheadrightarrow\mspec(R_K)$ is a surjection, and
			\item\label{point:3_prufer_domain} the ring $R$ is the increasing union of its subrings that are semilocal Pr\"ufer domain of finite Krull dimension. 
		\en
	\el
	\bs
		The first claim is a consequence of the fact that any integral domain is the intersection of its localisations at the maximal ideals.
		\par \eqref{point:2_prufer_domain}$\colon$ This is a consequence of \cite[Theorem 12.2]{matsumura_comm_ring_theory} (or \cite[Chapter~VI, \textsection7, No.~1, Proposition~2]{bourbaki}). Indeed, by \cite[\href{https://stacks.math.columbia.edu/tag/0AAV}{Tag 0AAV}]{stacks-project}, the intersection $R_{\fr{m}}\cap K$ is a valuation ring with fraction field $K$, for each maximal ideal $\fr{m}\subset R$. Therefore, \bud R_{K}\colonequals R\cap K=\textstyle\bigcap_{\fr{m}\in\mathrm{MaxSpec}(R)} (R_{\fr{m}}\cap K)\eud is a semilocal Krull domain with fraction field $K$, namely, there are a surjection $\mspec(R)\twoheadrightarrow\mspec(R_K)$ given by $\fr{m}\mapsto\fr{m}_K\colonequals\fr{m}\cap R_K$ and an isomorphism $(R_K)_{\fr{m}_K}\cong R_{\fr{m}}\cap K$ for each maximal ideal.
		\par \eqref{point:3_prufer_domain}$\colon$ The fraction field $K$ of $R$ can be written as $K=\bigcup K'$, the increasing union being taken over subfields $K'\subset K$ that are finitely generated extensions over the prime subfield $\bb{F}\subset K$; as a consequence, $R=\bigcup (R\cap K')$. Thanks to \eqref{point:2_prufer_domain}, the semilocal Pr\"ufer domain $R'\colonequals R\cap K'$ has fraction field $K'$; therefore, the fact that $\trdeg_{\bb{F}}(K')<\infty$ ensures that $R'$ has finite Krull dimension (\ccite{bourbaki}{Chapter~VI, Section~10, Number~3, Corollary~1}). Therefore, we have exhibited $R$ as an increasing union of subrings $R'$ that are semilocal Pr\"ufer domain of finite Krull dimension. Hence, we are done.
	\es
	By \Cref{lem:prufer_domain}, the notion of semilocal Krull domains coincides with the same of semilocal Pr\"ufer domains. Consequently, as a result of \Cref{lem:pemanence_properties}, the category of semilocal Krull domains is closed under localisations and under quotients by prime ideals.
	\bl\label{lem:Prufer_domain_gluing}
		Let $R_1$ be a semilocal Pr\"ufer domain with a maximal ideal $\frak{m}$ and a residue field $k$ and let $R_2$ be a semilocal Pr\"ufer domain whose fraction field is $k$. The subring \bud\label{diag:lem:Prufer_domain_gluing}\textstyle R\colonequals R_1\times_{k}R_2\subseteq R_1\eud is a semilocal Pr\"ufer domain whose residue field at $\frak{p}\colonequals\frak{m}\cap R$ is $k$. Moreover, the semilocalisation of $R$ at $\frak{p}$ and the maximal ideals not containing $\frak{p}$ is $R_1$ and $R/\frak{p}\cong R_2$.
	\el
	\bs
		We begin by establishing that $R$ is a semilocal Pr\"ufer domain. According to \Cref{lem:prufer_domain}, this is equivalent to demonstrating that it is a semilocal Krull domain. The same lemma shows that the rings $R_1$ and $R_2$ are semilocal Krull domains. To prove that $R$ itself is a semilocal Krull domain, we can, without loss of generality, reduce to the case where $R_1$ and $R_2$ are valuation rings. Indeed, this reduction is possible thanks to the fact that intersections of rings commute with products of rings. The claim when $R_1$ and $R_2$ are valuation rings follows from \ccite{fujiwara_kato_foundations_rigid_geometry}{Chapter~0, Proposition~6.4.1(1)} (cf.~\stacks{088z}).
		\par Therefore, it suffices to demonstrate that $\frak{p}\subset R$ satisfies the remaining claims. Since, by \stacks{0b7j}, we have that $\spec(R)=\spec(R_1)\sqcup_{\spec(k)}\spec(R_2)$, the claim that $R/\frak{p}\cong R_2$ follows. Furthermore, the same coproduct description of $\spec(R)$ demonstrates that the semilocalisation of $R$ at $\frak{p}$ and the maximal ideals not containing $\frak{p}$ is $R_1$. Thus, the proof is complete.
	\es
	\br
		Given a ring $R$, we may produce a graph with vertices in $\spec(R)$, where two vertices are joined by an edge if either is an immediate specialisation of the other. For example, the graph of a valuation ring $R$ is a path with two endpoints, which corresponds to the zero and the maximal ideal of $R$ respectively. Consequently, by the compatibility of forming this graph with localisation of rings, it is clear that the graph of a Pr\"ufer domain is a connected, acyclic graph. Conversely, an application of \Cref{lem:Prufer_domain_gluing} demonstrates that any finite, connected, acyclic graph may be recursively realised as the graph of a semilocal Pr\"ufer domain of finite Krull dimension.
	\er
	\bl\label{lem:Prufer_domain_cutting_and_pasting}
		Given a semilocal Pr\"ufer domain $R$ of finite Krull dimension, a prime ideal $\frak{p}\subset R$ and an element $a\in R$ such that $V(a)=\{\frak{q}\in\spec(R)\mid\frak{q}\supsetneq\frak{p}\}$, the canonical morphism \bd\label{diag:lem:Prufer_domain_cutting_)and_pasting}\textstyle R\arrow[r,"\sim"]&(R\inv{a})\times_{(R\inv{a}/(\frak{p}R\inv{a}))}(R/\frak{p})\text{\hspace{0.35cm}is an isomorphism.}\ed
	\el
	\begin{rems}\label{rems:prime_voidance}~
		\bn
			\item\label{part:rem:prime_avoidance} Since $R$ has finitely many prime ideals, such an $a$ always exist thanks to prime avoidance \stacks{00ds}.
			\item The element $a$ is chosen such that $R\inv{a}$ is the semilocalisation of $R$ at $\frak{p}$ and the maximal ideals not containing $\frak{p}$. In particular, $\frak{p}R\inv{a}\subset R\inv{a}$ is a maximal ideal.
		\en
	\end{rems}
	\bs[Proof of \Cref{lem:Prufer_domain_cutting_and_pasting}]
		The claim \eqref{diag:lem:Prufer_domain_cutting_)and_pasting} when $R$ is a valuation ring is the content of \ccite{fujiwara_kato_foundations_rigid_geometry}{Chapter~0, Proposition~6.4.1(1)}. Therefore, without loss of generality, we may assume that $R$ is not a valuation ring. Since $R$ is an integral domain, we have that $R\subseteq R\inv{a}$; consequently, the canonical morphism in claim \eqref{diag:lem:Prufer_domain_cutting_)and_pasting} is an injection. We need to show that this morphism is a bijection. Letting $R'\colonequals (R\inv{a})\times_{(\kappa(\frak{p}))}(R/\frak{p})$, by \stacks{0b7j}, in similar vein as the proof of \Cref{lem:Prufer_domain_gluing}, we have that \bud\textstyle\spec (R')=\spec(R\inv{a})\sqcup_{\spec(\kappa(\frak{p}))}\spec (R/\frak{p}).\eud However, the latter topological space is homeomorphic to $\spec(R)$. Therefore, the maximal ideals (in reality, all the prime ideals) of $R'$ and $R$ correspond to each other. Moreover, by \Cref{lem:pemanence_properties}, each multiplicative factor of $R'=(R\inv{a})\times_{(\kappa(\frak{p}))}(R/\frak{p})$ is a semilocal Pr\"ufer domain. As a result, by \Cref{lem:Prufer_domain_gluing}, we get that $R'$ itself is a semilocal Pr\"ufer domain. Since the maximal ideals of $R'$ and $R$ correspond and each ring is a semilocal Krull domain (\Cref{lem:prufer_domain}), to verify that the morphism in claim \eqref{diag:lem:Prufer_domain_cutting_)and_pasting} is an isomorphism, it suffices to ensure that their localisations at each maximal ideal are equal (\ccite{matsumura_comm_ring_theory}{Theorem~12.2} or \ccite{bourbaki}{Chpater~6, Section~7, Number~1, Proposition~2}). We verify this claim below. Let $\frak{m}\subseteq R$ be a maximal ideal. If $\frak{m}\not\supseteq\frak{p}$, we get that \bud\textstyle(R')_{(\frak{m}R')}=(R\inv{a})_{(\frak{m}R\inv{a})}=R_\frak{m}.\eud Therefore, it remains to verify the case when $\frak{m}\supseteq\frak{p}$. In this case, \bud\textstyle(R')_{(\frak{m}R')}=(R\inv{a})_{(\frak{m}R\inv{a})}\times_{(\kappa(\frak{p}))}(R/\frak{p})_{(\frak{m}/\frak{p})}=R_{\frak{p}}\times_{(\kappa(\frak{p}))}(R_{\frak{m}}/(\frak{p}R_{\frak{m}})).\eud Given that, by hypothesis, $R_\frak{m}$ is a valuation of finite rank (whose spectrum forms a finite linear chain), the second equality follows from the fact that $(R\inv{a})_{(\frak{m}R\inv{a})}=(R_{\frak{m}})\inv{a}=R_{\frak{p}}$. On the other hand, the valuation ring case of the claim \eqref{diag:lem:Prufer_domain_cutting_)and_pasting} produces the equality \bud R_{\frak{p}}\times_{(\kappa(\frak{p}))}(R_{\frak{m}}/(\frak{p}R_{\frak{m}}))=R_{\frak{m}}.\eud As a result, we get that $R_{\frak{m}}=(R')_{\frak{m}R'}$, as claimed. Thus, we are done.
	\es
	The following is an input in the proof of \Cref{lem:approximation:prufer_domain} below.
	\bl\label{lem:valuation_ring_pre-image}
		Let $R$ be an integral domain, let $\frak{m}\subset R$ be a maximal ideal with a residue field $k$, and let $\ell\subseteq k$ be a subfield. Suppose that $\varphi\colon R\twoheadrightarrow k$ is the canonical surjection. The subring \bud\text{$R_{\ell}\colonequals\{r\in R\mid \varphi(r)\in\ell\}$ is an integral subdomain such that $R_{\ell}=\Frac(R_{\ell})\cap R$.}\eud
	\el
	\bs
		Let $R'_{\ell}\colonequals\Frac(R_{\ell})\cap R$ and $\ell'\colonequals\varphi(R'_{\ell})$. To demonstrate that $R_{\ell}=R'_{\ell}$, it suffices to prove that $R'_{\ell}\subseteq R_{\ell}$. Applying $\varphi$, it is equivalent to show that $\ell'\subseteq\ell$. This follows since $\ell$ is a field and since $R'_{\ell}\subseteq\Frac(R_{\ell})$.
	\es
	We are now ready to prove a generalisation of \Cref{lem:prufer_domain}\eqref{point:3_prufer_domain}, which can be recovered by setting $\frak{p}=(0)$ in \Cref{lem:approximation:prufer_domain}\eqref{claim:ii}. We introduce some terminology below for convenience.
	\par Given a semilocal Pr\"ufer domain $R$, a set of prime ideals $\mathcal{P}=\{\frak{p}_r\subset R\}$ and a subring $R'\subseteq R$, we say that $R'$ is \textit{finitely generated at $\cal{P}$} if the residue field of $R'$ at $\frak{p}_r\cap R'$ is finitely generated over its prime subfield, for each $\frak{p}_r\in\cal{P}$. 
	\bl\label{lem:approximation:prufer_domain}
		Let $R$ be a semilocal Pr\"ufer domain of finite Krull dimension with a fraction field $K$ and let $\mathcal{P}\colonequals\{\frak{p}_r\subset R\}$ be a set of prime ideals.
		\bn[(i)] 
			\item\label{claim:i} Given a subfield $K'\subseteq K$, if $R$ is finitely generated at $\cal{P}$, then same is true for the subring $R'\colonequals R\cap K'$.
			\item\label{claim:ii} The ring $R$ is a filtered union of semilocal Pr\"ufer domains $R_{\alpha}$ of finite Krull dimension each of which is finitely generated at $\cal{P}$.
			\item\label{claim:ii'} The field $K$ can be written as a filtered union of subfields $K_{\alpha}\subseteq K$ so that each $R_{\alpha}\colonequals R\cap K_{\alpha}$ is finitely generated at $\cal{P}$. 
		\en
	\el
	\br
		\ccite{bourbaki}{Chapter~VI, \textsection~10, No.~3, Corollary~1} demonstrates that the transcendence degrees of $\kappa(\frak{p}_{r})$ over their respective prime subfields are bounded, for all $r$, as soon as $K$ is finitely generated over its prime subfield. However, ensuring that $\kappa(\frak{p})$ are also finitely generated is more delicate.
	\er
	\bs[Proof of \Cref{lem:approximation:prufer_domain}]
		First, we observe that \eqref{claim:ii} is implied by \eqref{claim:ii'}. Indeed, by \Cref{lem:prufer_domain}\eqref{point:2_prufer_domain}, such subrings $R_{\alpha}\colonequals R\cap K_{\alpha}$ are automatically semilocal Pr\"ufer domains of finite Krull dimension. Therefore, it suffices to demonstrate only \eqref{claim:i} and \eqref{claim:ii'}. We shall prove them by inducting on $r$. 
		\par The claims are trivial in the base case, i.e., when $r=0$. Therefore, we may assume that $r\ge 1$. As inductive hypothesis, we suppose that both \eqref{claim:i} and \eqref{claim:ii'} are true for all such $R$ when $r=n$. As inductive step, we shall prove both the claims when $r=n+1$. Let $\frak{p}\in\mathcal{P}$ be a maximal element ordered by inclusion and let $a\in R$ be such that $V(a)=\{\frak{q}\in\spec(R)\mid\frak{q}\supsetneq\frak{p}\}$. Since $R$ has only finitely many prime ideals, the element $a$ exists by prime avoidance \stacks{00ds} (cf.~\Cref{rems:prime_voidance}\eqref{part:rem:prime_avoidance}). Thanks to \Cref{lem:Prufer_domain_cutting_and_pasting}, \bd\label{diag:R_breakdown_at_p}\textstyle R\arrow[r,"\sim"]&(R\inv{a})\times_{(\kappa(\frak{p}))}(R/\frak{p})\text{\hspace{0.35cm}is an isomorphism.}\ed Furthermore, $R\inv{a}$ is a semilocal Pr\"ufer domain of finite Krull dimension whose residue field at the maximal ideal $\frak{p}R\inv{a}$ is $\kappa(\frak{p})$ and $R/\frak{p}$ is a semilocal Pr\"ufer domain of finite Krull dimension whose fraction field is $\kappa(\frak{p})$ (\Cref{lem:pemanence_properties}). Again, considering the fact that $R$ has only finitely many primes, it follows that $\spec(R_{\frak{p}})\subseteq\spec(R)$ is a Zariski open. As a consequence, since both the claims are Zariski local in $R$, by localising $R$ at $\frak{p}$, we may, without loss of generality, assume that $R$ is a valuation ring with a maximal ideal $\frak{p}\subset R$. 
		First, we prove \eqref{claim:i}.\vspace{0.3cm}
		\par\eqref{claim:i}$\colon$ Let $K'\subseteq K$ be subfield. As a first case, we assume that $\cal{P}=\{\frak{p}\}$. Since $R'\subseteq R$, the residue field $k'$ of $R'$ at $\frak{p}\cap R'$ is a subfield of $\kappa(\frak{p})$, which is finitely generated over its prime subfield. Consequently, by \ccite{bourbaki_algebraII_Chapters4-7}{Chapter~V, \textsection14, No.~7, Corollary~3}, the same is true for $k'$. Thus, in this case, we are done.
		\par Otherwise, the set $\cal{P}'\colonequals(\cal{P}\setminus\{\frak{p}\})$ is nonempty. Let $\frak{q}$ be the maximal element in $\cal{P}'$ (the prime $\frak{q}$ is unique because $\spec(R)$ is totally ordered). Similar to \eqref{diag:R_breakdown_at_p}, we have an isomorphism \bd\label{diag:R_breakdown_at_q}\textstyle R\arrow[r,"\sim", "\eqref{diag:lem:Prufer_domain_cutting_)and_pasting}"']&R_{\frak{q}}\times_{(\kappa(\frak{q}))}(R/\frak{q}).\ed Letting $\frak{q}'\colonequals\frak{q}\cap R'$, since $R'_{\frak{q}'}$ is a valuation ring that is dominated by $R_{\frak{q}}\cap K'$, we obtain the equality $R'_{\frak{q}'}=R_{\frak{q}}\cap K'$. Let $k'$ be the residue field of $R'$ at $\frak{q}'$. Since $R$ is finitely generated at $\cal{P}$, it follows that $R_{\frak{q}}$ (respectively, $R/\frak{q}$) is finitely generated at $\cal{P}'$ (respectively, at $\{(\frak{p}/\frak{q})\}$). Consequently, by the induction hypothesis (respectively, by the $r=1$ case proved above), the ring $R'_{\frak{q}'}$ (respectively, the ring $R'/\frak{q}'$) is finitely generated at $\cal{P}'$ (respectively, at $\{(\frak{p}/\frak{q})\}$). Therefore, the ring \bud\textstyle R'\arrow[r,"\sim", "\eqref{diag:lem:Prufer_domain_cutting_)and_pasting}"']&R'_{\frak{q}'}\times_{k'}(R'/\frak{q}')\eud is finitely generated at $\cal{P}$, as required. The induction is thus complete and \eqref{claim:i} is proven.
		\vspace{0.3cm}\newline Finally, we prove \eqref{claim:ii'}.\vspace{0.3cm}
		\par\eqref{claim:ii'}$\colon$ In a similar vein as the proof of \eqref{claim:i}, as a first case, we assume that $\cal{P}=\{\frak{p}\}$. Letting $\bb{F}\subseteq\kappa(\frak{p})$ be the prime subfield, we can write $\kappa({\frak{p}})$ as a filtered unions of subfields $k_{\beta}\subseteq\kappa(\frak{p})$ which are finitely generated over $\bb{F}$. For each $\beta$, we let \bud\textstyle R_{\beta}\colonequals\{r\in R\mid r~\mathrm{mod}~\frak{p}\in k_{\beta}\}\hspace{0.5cm}\text{ and }\hspace{0.5cm}K_{\beta}\colonequals\Frac(R_{\beta}).\eud By Lemmas~\ref{lem:valuation_ring_pre-image} and \ref{lem:prufer_domain}\eqref{point:3_prufer_domain}, the ring $R_{\beta}$ is a semilocal Pr\"ufer domain of finite Krull dimension whose residue field at the maximal ideal $\frak{m}_{\beta}\colonequals \frak{p}\cap R_{\beta}$ is $k_{\beta}$, for each $\beta$. Furthermore, the same lemmas prove the equality $R_{\beta}=R\cap K_{\beta}$, for each $\beta$. As a consequence, it suffices to check that $\bigcup K_{\beta}=K$, which is true because $\bigcup k_{\beta}=k$. Thus, in this case, we are done.
		\par Otherwise, the set $\cal{P}'\colonequals \cal{P}\setminus\{\frak{p}\}$ is nonempty. Similar to the proof of \eqref{claim:i}, let $\frak{q}$ be the maximal element in $\cal{P}'$. By the induction hypothesis (respectively, by the $r=1$ case proved above), the field $K$ can be written as a filtered union of subfields $K_{\gamma}\subseteq K$ (respectively, subfields $K_{\beta}\subseteq K$) such that $R_{\gamma}\colonequals R\cap K_{\gamma}$ (respectively, $R_{\beta}\colonequals R\cap K_{\beta}$) is finitely generated at $\cal{P}'$ (respectively, at $\{\frak{p}\}$), for all $\gamma$ (respectively, for all $\beta$). Possibly by reindexing, we consider a mono-indexed cofinal collection $\{K_{\alpha}\}$ in the bi-indexed collection $\{K_{\beta}\cap K_{\gamma}\}$ of subfields of $K$ and we let \bud\text{$R_{\alpha}\colonequals R\cap K_{\alpha}$, for each $\alpha$}.\eud Since $\{K_{\alpha}\}$ is a cofinal collection of $\{K_{\beta}\cap K_{\gamma}\}$, letting $\alpha$ be an index, there exists some $\beta$ and some $\gamma$ (we fix one such pair) so that $K_{\alpha}=K_{\beta}\cap K_{\gamma}$. Therefore, by definition, we get \bud R_{\alpha}=R_{\beta}\cap K_{\gamma}=R_{\gamma}\cap K_{\beta}.\eud Applying \eqref{claim:i} to $R_{\gamma}$ (respectively, to $R_{\beta}$), we obtain that $R_{\alpha}$ is finitely generated at $\cal{P}'$ (respectively, at $\{\frak{p}\}$). Combining the statements, it follows that $R_{\alpha}$ is finitely generated at $\cal{P}=\cal{P}'\cup\{\frak{p}\}$, as required. 
		Thus, the induction step is complete, and consequently, the claim is proven. Hence, we are done.
	\es
	Below, in \Cref{lem:snail_generic_point_of_special_fibre}\eqref{claim:1}, we recall from \ccite{phd-thesis}{Lemma~3.10} (cf.~\ccite{mb_valuation_ring}{théorème~A}) an intriguing property of smooth algebras $A$ over valuation rings $V$. Specifically, we demonstrate that the local rings of $A$ at the generic points of its $R$-fibres are also valuation rings. We will leverage this to ultimately reduce \Cref{thm:1.1} for $A$ to a Gersten's injectivity claim for the local ring of $A$ at a generic point of its $R$-special fibre (see \Cref{cor:Gersten_injectivity_for_valuation}).
	\par Before stating \Cref{lem:snail_generic_point_of_special_fibre}, we introduce some notations for clarity. Given a ring $R$, we denote its Krull dimension by $\dim(R)$. Given a finite set of primes $\ca{P}\subseteq\spec(A)$ of a ring $A$, we denote the semilocalisation of $A$ at the primes in $\cal{P}$ by $A_{\cal{P}}$. Given an extension $e\colon R\to S$ of rings and a subset $\ca{P}\subset\spec(S)$, let \bud\mathrm{Flat}_{\mathcal{P}}(S/R)\colonequals\{f\in S\mid V(f)\cap\mathcal{P}=\emptyset \text{ and } R\to S/fS\text{ is flat}\}.\eud
	In terms of the previously introduced notations, for any extension $e$ and subset $\mathcal{P}$, it is noteworthy that $\mathrm{Flat}_{\mathcal{P}}(S/R)$ is a filtered set. 
	\bl\label{lem:snail_generic_point_of_special_fibre}~
		\bn[(1)]
			\item\label{claim:1} Given a valuation ring $V$, an integral domain $A$ that is a smooth $V$-algebra, and a subset $\mathcal{P}$ of the set $\tilde{\mathcal{P}}$ of generic points of the $V$-special fibre of $\spec(A)$, the morphism $V\to A_{\cal{P}}$ is faithfully flat extension of semilocal Pr\"ufer domains.
			\item\label{claim:2} Moreover, if $\dim(V)<\infty$, then there exists an affine open neighbourhood $\spec(A')\subseteq\spec(A)$ of $\cal{P}$ such that we have \bud\colim_{f\in\mathrm{Flat}_{\mathcal{P}}(A'/V)} A'\inv{f}=A'_{\ca{P}}.\eud
		\en
	\el
	\br\label{rem:generic_fibre}~
		\bn
			\item Although, our lemma is stated for a general $\cal{P}$, the most interesting cases are when $\cal{P}$ is a singleton set and when $\cal{P}=\tilde{\cal{P}}$. In the former case, the proof of \eqref{claim:1} below shows that, in fact, $V\to A_{\cal{P}}$ is an extension of valuation rings which induces an isomorphism between the respective value groups.
			\item\label{point:remark:generic_fibre} It is worth noting that the $R$-special fibre of $\spec(A)$ need not be connected, even when $R$ is a discrete valuation ring and $R\to\cal{A}$ is \'etale. For example, when $\cal{A}$ is the integral closure of $R$ is an unramified extension $L\supseteq\Frac(R)$ (\stacks{09e9}), the maximal ideals of $\cal{A}$, which lie in the $R$-special fibre of $\spec(A)$, correspond to the extensions of valuations $R\subset L$ centred on $\cal{A}$ (\stacks{09e8}). In particular, if $R$ is not Henselian, such extensions of valuations need not be unique.  
		\en
	\er
	\bs[Proof of \Cref{lem:snail_generic_point_of_special_fibre}]
		We prove the claims simultaneously. We shall reduce, by the local structure of smooth morphisms, to establishing the claims when $A$ is \'etale, and when $A$ is a polynomial algebra, in which we do an explicit computation. 
		\par By the semilocal structure of smooth morphisms \ccite{sga1}{Expos\'e II, th\'eor\`eme 4.10(ii)}\footnote{Even though \ccite{sga1}{Expos\'e II, th\'eor\`eme 4.10(ii)} concerns only the Zariski local structure of smooth morphisms, we can establish their Zariski semilocal structure with a slight modification of the proof in loc.~cit.} (cf.~\cite[\href{https://stacks.math.columbia.edu/tag/052E}{Tag 052E}]{stacks-project}), there exist an integer $n\ge 0$, a connected, affine, open neighbourhood $\spec(B)\subseteq\spec(A)$ of $\cal{P}$ and an \'etale morphism $j\colon\spec(B)\to\spec(V[x_1,\ldots,x_n])$. Since the claims are of Zariski semilocal in nature, without loss of generality, we may assume that $A=B$. If $V\to A$ is \'etale, then $n=0$, and, by \cite[\href{https://stacks.math.columbia.edu/tag/0ASJ}{Tag 0ASJ}]{stacks-project}, claim \eqref{claim:1} is true. Furthermore, if $\dim(V)<\infty$, the ring $A$ is a semilocal Pr\"ufer domain of finite Krull dimension. Therefore, by prime avoidance \stacks{00ds}, $\spec(A_{\ca{P}})\subseteq\spec(A)$ is an open subset. Consequently, to demonstrate claim \eqref{claim:2}, it suffices to take $A'=A_{\cal{P}}$. Thus, the claims are proven in the case when $n=0$.\vspace{0.25cm} 
		\par Therefore, without loss of generality, we may assume that $n\ge 1$. Since \'etale morphisms are of relative dimension $0$, it follows that $j(\frak{p})$ is the generic point $\eta$ of the $V$-special fibre of $\spec(V[x_1,\ldots,x_n])$, for any $\frak{p}\in\cal{P}$. As a consequence, letting $V'\colonequals V[x_1,\ldots,x_n]$ and $\frak{p}'\subset V'$ be the prime corresponding to $\eta$, we have a factorisation \bud V\xrightarrow{g} V'_{\fr{p}'}\xrightarrow{h} A_{\ca{P}}.\eud Thus, it is sufficient to show the claims are true for the morphism $g$. Indeed, supposing that the claims are true for $g$, the proven $n=0$ case applied to $h$ implies that $h\circ g$ is a faithfully flat extension of semilocal Pr\"ufer domains, as required to show \eqref{claim:1}. Similarly, to prove claim \eqref{claim:2}, we note that there is an inclusion $\mathrm{Flat}_{\{\fr{p}\}}(V'/V)\subseteq\mathrm{Flat}_{\ca{P}}(A/V)$ and base change the equality obtained by \eqref{claim:2} for $g$ along $V'\to A$ to obtain an equality \bd\label{diag:case_2}\textstyle A_{\frak{p}'}=\colim_{f\in\mathrm{Flat}_{\{\fr{p}'\}}(V'/V)}A\inv{f}=\colim_{f\in\mathrm{Flat}_{\ca{P}}(A/V)}A\inv{f}.\ed Finally, base changing the equality \eqref{diag:case_2} along the open embedding $\spec(A_{\ca{P}})\subset\spec(A)$, we obtain the required equality which shows \eqref{claim:2}. Therefore, without loss of generality, it suffices to establish the claims for $g$. 
		\par In consequence, we may assume that $A=V[x_1,\ldots,x_n]$. Let $\fr{m}\subset V$ be the maximal ideal and let $\fr{p}\colonequals \fr{m}[x_1,\ldots,x_n]$. Given the assumption, we have that $\ca{P}=\{\fr{p}\}$. \vspace{0.25cm}
		\par\eqref{claim:1}$\colon$ In fact, we shall demonstrate that $V\to A_{\fr{p}}$ is a faithfully flat extension of valuation rings that induces an isomorphism an isomorphism between the respective value groups. First, we prove that $A_{\fr{p}}$ is a valuation ring. To do so, it suffices to verify that for any \bud\text{$t\in\Frac(A_{\fr{p}})=\Frac(V[x_1,\ldots,x_n])$, either $t\in A_{\fr{p}}$ or $1/t\in A_{\fr{p}}$.}\eud Let $t=f/g\in\Frac(A_{\fr{p}})$, where $f,g\in A$. Using the valuation on $V$, we define $\mathrm{val}(f)\in V$ (resp., $\mathrm{val}(g)\in V$) to be the element, which is well defined up to a unit in $V$, such that $f/\mathrm{val}(f)\in A\setminus\fr{p}$ (resp., $g/\mathrm{val}(g)\in A\setminus\fr{p}$). If $\mathrm{val}(g)\mid\mathrm{val}(f)$, then $t\in A_{\fr{p}}$, otherwise, $1/t\in A_{\fr{p}}$, and we are done. It remains to show that the morphism \bud\text{$\varphi\colon\units{\Frac(V)}/\units{V}\to\units{\Frac(A)}/\units{A_{\fr{p}}}$ of value groups is an isomorphism.}\eud Since $\varphi$ is injective, it suffices to show that $\varphi$ is surjective. In a similar vein as the previous arguments, given $t\in\units{\Frac(A)}$, there exists $u\in\units{\Frac(V)}$ such that $t/u\in\units{A_{\fr{p}}}$. Thus, \eqref{claim:1} is proven.\vspace{0.25cm}
		\par\eqref{claim:2}$\colon$ Since, by definition, $A_{\fr{p}}=\colim_{f\in A\setminus\{\fr{p}\}}A\inv{f}$, it is enough to show that $A\setminus\{\fr{p}\}\subseteq\mathrm{Flat}_{\{\fr{p}\}}(A/V)$. In other words, given $f\in A\setminus\{\fr{p}\}$, we need to demonstrate that $V\to A/fA$ is flat. As $V$ is a valuation ring, it is equivalent to prove that such an $f$ satisfies that $A/fA$ is $V$-torsion free. Letting $f\in A\setminus\{\fr{p}\}$, suppose that $v\in V$ and $a,b\in A$ be any elements such that $va=bf$. In order to establish that $A/fA$ is $V$-torsion free, we need to show that either $v=0$ or $a\in fA$. Therefore, further assuming that $v\neq 0$, it is sufficient to prove that $a\in fA$. Adopting a notation of the above paragraph and keeping in mind that $f\not\in\fr{p}$, we conclude that $v$ divides $\mathrm{val}(bf)=\mathrm{val}(b)$. As a consequence, $b'\colonequals b/v\in A$, whence we deduce that $a=b'f$, as required. Thus, \eqref{claim:2} is also proven, and hence, we are done. 
	\es
	\section{Gersten's Injectivity of $G$-theory}\label{section:G-theory}
	In this section, our goal is to prove \Cref{thm:4}. Although the anatomy of its proof is based on the proof of \ccite{gillet_levine}{Theorem}, considerable amount of customisation is needed to upgrade the arguments to the non-Noetherian setting. Already, since the spectra of valuation rings contain more points than their discrete counterparts, the construction of Presentation Lemma~\ref{lem:presentation_lemma_Prufer}, which is the principal technical ingredient in this section, is surprisingly more delicate than we expect (cf.~\ccite{gillet_levine}{Lemma~1} and \ccite{luders_relative_gersten_conjecture_for_milnor_k_theory}{Lemma~2.12}). We assemble our version of the presentation lemma by synthesising the proofs of \ccite{ces_grothendieck-serre}{Variant~3.7} and \ccite{phd-thesis}{Proposition~6.4}. In conclusion of this section, and as an application of \Cref{thm:1.3}, we extend the result in \ccite{gillet_levine}{Theorem} to demonstrate that a functor $\scr{F}$, which satisfies the criteria outlined in \Cref{defn:G-theory_like_functor}—namely, it commutes with filtered colimits of rings, is a localizing invariant, and has pushforwards along finite morphisms—satisfies Gersten's injectivity for smooth algebras over equicharacteristic valuation rings (\Cref{cor:Gersten_injectivity_for_valuation}). 
	\par In preparation for Presentation Lemma~\ref{lem:presentation_lemma_Prufer}, let us state a straightforward consequence of \ccite{ces_grothendieck-serre}{Proposition~3.6}. 
	\bl\label{lem:presentation_lemma_field}
		Given a field $k$, an affine, smooth $k$-scheme $X$ of pure relative dimension $d>0$, points $x_1,\ldots, x_n\in X$ and a nowhere dense, closed subscheme $Y\hookrightarrow X$, there exist affine opens $x_1,\ldots,x_n\in U\subseteq X$ and $S\subseteq\bb{A}^{d-1}_k$ and a smooth $k$-morphism $\pi\colon U\to S$ of relative dimension $1$ such that $\pi|_{Y\cap U}$ is finite.
	\el
	\bs
		We shall deduce the claim via a direct application of loc.~cit. At the cost of shrinking $X$ to an open, affine neighbourhood of $x_1,\ldots,x_n$, we may assume that there is a closed embedding $X\hookrightarrow\bb{A}^m_k$, and thus, by composing with a chosen open embedding $\bb{A}^m_k\hookrightarrow\bb{P}^m_k$, we obtain an embedding $\iota\colon X\hookrightarrow\bb{P}^m_k$. We define $\overline{X}$ (respectively, $\overline{Y}$) to be the schematic closure of $\iota$ (respectively, the same of $Y\subset\overline{X}$). The fact that $Y$ is dense in $\overline{Y}$ ensures that $\codim_{\overline{X}}(\overline{Y})=\codim_{X}(Y)\ge 1$. Moreover, the same fact implies that $\overline{Y}\setminus Y$ does not contain any generic point of $\overline{Y}$. As a consequence, it follows that $\codim_{\overline{X}}(\overline{Y}\setminus Y)\ge 2$. Therefore, the claim follows by putting ($\overline{X}, X, x_1,\ldots, x_n\in X, \overline{Y}$) in loc.~cit.
	\es
	\begin{pl}\label{lem:presentation_lemma_Prufer}
		Given 
		\bun
			\item[$\circ$] a semilocal Pr\"ufer domain $R$ of finite Krull dimension, 
			\item[$\circ$] an affine, smooth $R$-scheme $X$ of pure relative dimension $d>0$, 
			\item[$\circ$] points $x_1,\ldots, x_n\in X$ and 
			\item[$\circ$] an $R$-flat, closed subscheme $Y\hookrightarrow X$ that does not contain any component of the $R$-fibres of $X$, 
		\eun
		there exist affine opens $x_1,\ldots,x_n\in U\subseteq X$ and $S\subseteq\bb{A}^{d-1}_R$ and a smooth $R$-morphism $\pi\colon U\to S$ of relative dimension $1$ such that $\pi|_{Y\cap U}$ is quasi-finite.
	\end{pl}
	For brevity, such a morphism $\pi$ (along with the data of $U$ and $S$) will be called a \textit{presentation of $X$ over $R$.}
	\bs
		Let $C\subseteq\spec R$ be the closed subscheme of closed points and let $X_C$ be the fibre of $X$ over $C$. Since \Cref{lem:presentation_lemma_field} takes care of the case when $R$ is a field, without loss of generality, we may assume that $C\subsetneq\spec R$. \vspace{0.25cm}
		
		\par \underline{Case 1}: We suppose further that each $x_i$ specialises to a point in $X_C$. Because of this assumption, we may replace $x_i$ with one of its specialisations that are closed points in $X_C$ (we use the fact that $X_C$ is a Jacobson scheme in order to specialise each $x_i$ to a closed point in $X_C$) to assume further that each $x_i\in X_C$ is closed. Thanks to \Cref{lem:presentation_lemma_field}, there exist affine opens \bud\text{$x_1,\ldots,x_n\in U'\subset X_C$ and $S'\subset\bb{A}^{d-1}_C$ and a smooth $C$-morphism $\pi'\colon U'\to S'$}\eud such that $\pi'|_{Y_C\cap U'}$ is quasi-finite (in fact, we can ensure that the latter is finite, but we will not require it for the proof of the claim). Let $\scr{I}$ be the ideal of vanishing of $X_C\hookrightarrow X$ and let $\tilde{U}\subset X$ be an open subset such that $\tilde{U}\cap X_C=U'$. By lifting sections along $\scr{O}_{\tilde{U}}\twoheadrightarrow\scr{O}_{\tilde{U}}/(\scr{I}|_{\tilde{U}})$, we define a morphism $\tilde{\pi}\colon\tilde{U}\to\bb{A}^{d-1}_R$. Since $\tilde{U}$ is $R$-finite type and $R$-flat, it is $R$-finitely presented (by \ccite{raynaud_gruson}{Premi\`ere partie, théorème~3.4.6} and by the fact that $R$ is an integral domain). Thanks to the fibrewise criterion of flatness \stacks{039c}, $\tilde{\pi}$ is flat at each $x_i$, whence thanks to the fibrewise criterion of smoothness \stacks{01v8}, $\tilde{\pi}$ is smooth at each $x_i$. Therefore, after shrinking $\tilde{U}$, we may assume that $\tilde{\pi}$ is smooth. The openness of the quasi-finite locus \stacks{o1t1} implies that there exists an open subset $U_1\subset Y$ containing $Y\cap\pi^{-1}(\pi(x_i))$, for all $i$, such that $\pi|_{U_1}$ is quasi-finite. We choose affine opens $\pi(x_1),\ldots,\pi(x_n)\in S\subseteq\pi(\tilde{U}\setminus(Y\setminus U_1))$ and $U\subseteq\pi^{-1}(S)$.
		\vspace{0.25cm}
		\par \underline{Case 2}: In general, some of the points among $x_1,\ldots,x_n$ might not specialise to a point of $X_C$, say $y_1,\ldots,y_m$. Let $\mathcal{P}\subseteq\spec(R)$ be the images of $y_1,\ldots,y_m$. We shall tailor $X$ in such a way that each $y_i$ specialises to $X_C$, effectively replacing the original $X$ with this customised version. For this purpose, we use an ingenious trick due to \citeauthor{ces_grothendieck-serre} in the proof of \ccite{ces_grothendieck-serre}{Variant~3.7}. In this regard, thanks to \Cref{lem:approximation:prufer_domain}\eqref{claim:ii} and a limit argument, without loss of generality, we may assume that the residue field $\kappa(\frak{p})$ of $R$ at each $\frak{p}\in\cal{P}$ is finitely generated over its prime subfield. By, for example, \ccite{phd-thesis}{Lemma~6.1}, each field $\kappa(\frak{p})$ is a fraction field of a regular domain $A_{\frak{p}}$ that is smooth over $\bb{F}_p$ or $\bb{Z}$. Moreover, each $A_{\fr{p}}$ is of positive Krull dimension, since otherwise $K$ is a finite field, in which case, it contradicts our assumption that $R$ is not a field. By localising $A_{\fr{p}}$ and possibly by selecting a local $R$-projective embedding of $X$, we may assume that
		\bn
			\item the scheme $X_{\kappa(\frak{p})}$ spreads out to a smooth $A_{\frak{p}}$-scheme $X_{\frak{p}}$ that is fibrewise of pure dimension $d$ (see \stacks{01v8} and \ccite{egaIV_3}{th\'eor\`eme~12.1.1(iv)}),
			\item each point $y_i$ lying over $\frak{p}$ spreads out to an $A_{\frak{p}}$-finite, closed subscheme in $X_{\frak{p}}$, and
			\item the closed subscheme $Y_{\kappa(\frak{p})}$ spreads out to an $A_{\frak{p}}$-flat, closed subscheme $Y_{\frak{p}}$ such that $Y_{\frak{p}}$ is $A_{\frak{p}}$-fibrewise of codimension $\ge 1$ in $X_{\frak{p}}$ (see \stacks{039c} and \ccite{egaIV_3}{th\'eor\`eme~12.1.1(v)}).
		\en 
		Now that we have constructed the objects above for each $\fr{p}\in\ca{P}$, let us proceed with the customization of $X$.
		\vspace{0.25cm}
	 	\par \begin{whatever}{Method}{}\label{Method} Simultaneously,
	 		\bun 
	 			\item[$\circ$] we shall iteratively glue $R$ with a discrete valuation ring $A'_{\frak{p}}$ whose fraction field is $\kappa(\frak{p})$ at each $\frak{p}\in\cal{P}$ to ultimately produce a semilocal Pr\"ufer domain $\tilde{R}$ of finite Krull dimension $\tilde{R}$, and
	 			\item[$\circ$] we shall iterative glue $X$ with an $A'_{\frak{p}}$-scheme $X_{\frak{p}}$ whose generic fibre is $X_{\kappa(\frak{p})}$ at each $\fr{p}\in\cal{P}$ to ultimately produce an $\tilde{R}$-scheme $\tilde{X}$ which has the property that each $y_i$ specialises to a point of the $\tilde{R}$-special fibre of $\tilde{X}$.
	 		\eun
	 	\end{whatever}
	 	\vspace{0.25cm} Once this is done, we apply the already proven Case 1 to $\tilde{R}$ and $\tilde{X}$ (Case 1 applies thanks to our construction of $\tilde{R}$ and $\tilde{X}$) to obtain a presentation $\tilde{\pi}$ of $\tilde{X}$ over $\tilde{R}$. Given that a presentation of $X$ over $R$ is Zariski semilocal around $x_1,\ldots,x_n$, considering the open embedding $j\colon \spec(R)\hookrightarrow\spec(\tilde{R})$, we can perform a base change of $\tilde{\pi}$ along $j$ to obtain a presentation $\pi$ of $X$ over $R$, as required.
	 	\par Consequently, it suffices to construct such $\tilde{R}$ and such $\tilde{X}$. We proceed iteratively, incrementally considering primes $\frak{p}\in\cal{P}$ ordered by their height. We fix a prime $\frak{p}\in\cal{P}$ of height $n$. By abuse of notation, let $R$ be the semilocal Pr\"ufer domain (respectively, the $R$-scheme $X$) constructed by applying the procedure called Method to all $\frak{q}\in\cal{P}$ of height $<n$. Following Method, we shall construct $\tilde{R}$ by gluing $R$ with a discrete valuation ring $A'_{\frak{p}}$ at $\frak{p}$ and $\tilde{X}$ by gluing $X$ with an $A'_{\frak{p}}$-scheme $X_{\frak{p}}$ whose generic fibre is $X_{\kappa(\frak{p})}$. Let $a\in R$ be a element whose vanishing set is $\{\frak{q}\in\spec(R)\mid\frak{q}\supsetneq\frak{p}\}$. As observed in \Cref{lem:Prufer_domain_cutting_and_pasting}, thanks to prime avoidance \stacks{00ds}, such an $a$ always exists because $R$ has finitely many prime ideals. By \eqref{diag:lem:Prufer_domain_cutting_)and_pasting}, we may write \bud\textstyle R\arrow[r,"\sim"]&R\inv{a}\times_{\kappa(\frak{p})}(R/\frak{p}).\eud Given that $A_{\frak{p}}$ is of positive Krull dimension, it has infinitely many primes of height $1$, permitting us to choose such a prime $\frak{r}\subset A$ so that the localisation $A'_{\frak{p}}$, which is necessarily a discrete valuation ring, of $A_{\frak{p}}$ at $\frak{r}$ is different from each of the localisations of $R/\frak{p}$. Choosing such a prime $\frak{r}\subset A$, we substitute $A_{\frak{p}}$ with $A'_{\fr{p}}$ and consider $R' \colonequals (R/\frak{p})\cap A_{\frak{p}}$, where the intersection is taken in $\kappa(\frak{p})$. Thanks to \Cref{lem:prufer_domain}, the ring $R'$ is a semilocal Pr\"ufer domain of finite Krull dimension with fraction field $\kappa(\frak{p})$. We define $\tilde{R}$ by the following diagram \bd\label{diag:breakdown_definition_tilde_R}\textstyle \tilde{R}\arrow[r,"\sim"]&R\inv{a}\times_{\kappa(\frak{p})}R'.\ed \Cref{lem:Prufer_domain_gluing} ensures that $\tilde{R}$ is a semilocal Pr\"ufer domain of finite Krull dimension. It remains to construct $\tilde{X}$, which we do below. Over the open cover $\spec(R/\frak{p})$ and $\spec(A_{\frak{p}})$ of $\spec(R')$
		\bun 
			\item[$\circ$] we glue $X_{R/\frak{p}}$ and $X_{\frak{p}}$ along $X_{\kappa(\frak{p})}$ to obtain a smooth $R'$-scheme $X'$ that is fibrewise of pure dimension $d$, and 
			\item[$\circ$] we glue $Y_{R/\frak{p}}$ and $Y_{\frak{p}}$ along $Y_{\kappa(\frak{p})}$ to obtain an $R'$-flat, closed subscheme $Y'\subset X'$ which is $R'$-fibrewise of codimension $\ge 1$.
		\eun
		By construction, each $y_i$ lying over $\frak{p}$ specialises to a point in an $R'$-special fibre of $X'$. Finally, we define $\tilde{X}$. Thanks to \stacks{0b7j} and \stacks{0d2i}(1),
		\bun
			\item[$\circ$] we glue $X_{R\inv{a}}$ and $X'$ along $X_{\kappa(\frak{p})}$ to obtain an $\tilde{R}$-flat scheme $\tilde{X}$ that is fibrewise of pure dimension $d$, and
			\item[$\circ$] we glue $Y_{R\inv{a}}$ and $Y'$ along $Y_{\kappa(\frak{p})}$ to obtain an $\tilde{R}$-flat, closed subscheme $\tilde{Y}\subset\tilde{X}$ which is $\tilde{R}$-fibrewise of codimension $\ge 1$.
		\eun
		 Thanks to the fibrewise criterion of smoothness \stacks{01v8}, it follows that $\tilde{X}$ is $\tilde{R}$-smooth. By construction, the property that each $y_i$ lying over $\frak{p}$ specialises to a point in an $\tilde{R}$-special fibre of $\tilde{X}$ continues to be true. Therefore, we may apply Case 1 to ($\tilde{R}$; $\tilde{X}$; $x_1,\ldots,x_m$; $\tilde{Y}$) to obtain a presentation of $\tilde{X}$ over $\tilde{R}$, as required. This concludes the proof.
	\es
	\begin{rems}~
		\bn[(1)]
			\item Even though for the sake of demonstrating \Cref{thm:4}, we are only interested in a special case of Presentation Lemma~\ref{lem:presentation_lemma_Prufer}, namely, when $R$ is a valuation ring, the reason to state Presentation Lemma~\ref{lem:presentation_lemma_Prufer} in its current generality is not merely curiosity. As a matter of fact, our proof, more precisely, Case 2 of our proof forces us to modify $R$ via cute-and-glue tailoring (Lemmas~\ref{lem:Prufer_domain_gluing}-\ref{lem:Prufer_domain_cutting_and_pasting}). As a result, Pr\"ufer domains naturally appear there during this process. 
			\item There are significant differences between the statements of \ccite{phd-thesis}{Proposition~6.4} and Presentation Lemma~\ref{lem:presentation_lemma_Prufer}. Although our assumption that $Y$ is $R$-flat, as opposed to a fibrewise smoothness condition imposed on $Y$ in loc.~cit., makes our proof easier than there, a drawback is that we can only ensure that $\pi|_{Y\cap U}$ is quasi-finite, instead of finite. Quasi-finiteness is sufficient for the proof of \Cref{thm:4} thanks to the `linearity of $G$-theory', which greatly simplifies the arguments. Indeed, following the clever idea in the proof of \ccite{gillet_levine}{Theorem}, after base changing along $\pi$ (cf.~\eqref{diag:small_one}), we employ Zariski's main theorem \stacks{00q9} to enlarge our scheme, so that essentially $\pi|_{Y\cap U}$ may be assumed to be finite (cf.~\eqref{diag:big_one}), at least in some vague sense. 
		\en
	\end{rems}
	While working with non-Noetherian rings, we constantly need to make the distinction between finite type objects and finitely presented ones over them. This is vital, for example, especially while using Noetherian approximation techniques, where we require finitely presented, and not just finite type. Coherent rings, which we introduce below, form an important class of rings where this disparity between finite type and finite presentation reduces. Thankfully, rings in $\mathcal{S}_{\mathrm{val}}$ are coherent, which permits us to adapt some of the techniques from \cite{gillet_levine}.  
	\begin{montobo}[Coherence]\label{montobo:locally_coherent}
		Given a scheme $X$, an $\ca{O}_X$-module $\scr{F}$ is called \textit{coherent} if it is of finite type and for every open $U\subseteq X$ and every finite collection $s_i\in\scr{F}(U)$, $i=1,\ldots,n$, the kernel of the associated morphism $\bigoplus_{i=1,\ldots,n}\ca{O}_U\to\scr{F}$ is of finite type (\cite[\href{https://stacks.math.columbia.edu/tag/01BV}{Tag 01BV}]{stacks-project}). A coherent $\ca{O}_X$-module is finitely presented, and therefore, quasi-coherent (\cite[\href{https://stacks.math.columbia.edu/tag/01BW}{Tag 01BW}]{stacks-project}). A scheme $X$ is called \textit{locally coherent} if $\ca{O}_X$ is a coherent module over itself (\ccite{gabber-ramero-foundations}{Definition~8.1.54}). A ring $A$ is called \textit{coherent} if any finitely generated ideal of $A$ is finitely presented (\cite[\href{https://stacks.math.columbia.edu/tag/05CV}{Tag 05CV}]{stacks-project}). 
		\par A scheme $X$ that is locally of finite presentation over a Pr\"ufer domain $R$ is locally coherent. Indeed, since the property of being locally coherent is Zariski local, it suffices to check that any ring $A$ that is a finitely presented $R$-algebra is coherent. Let $f\colon A'\colonequals R[x_1,\ldots,x_n]\twoheadrightarrow A$ be a presentation of $A$ such that $\ker(f) \subset A'$ is a finitely generated ideal. Since $\ker(f) \subset A'$ is finitely generated, it is enough to show that the ring $A'$ is coherent. Letting $I\subset A'$ be a finitely generated ideal, we shall show that $I$ is a finitely presented $A'$-module. Putting $X=\spec A', S=\spec R$ and $\scr{M}=\widetilde{I}$ in \ccite{raynaud_gruson}{Premi\`ere partie, théorème~3.4.6} (by \ccite{bourbaki}{Chapter~I, \textsection2.4, Proposition~3(ii)}, \cite[\href{https://stacks.math.columbia.edu/tag/090Q}{Tag 090Q}]{stacks-project} and the fact that flatness is a local property \cite[\href{https://stacks.math.columbia.edu/tag/0250}{Tag 0250}]{stacks-project}, the $R$-torsion-free module $I$ is flat), we obtain that $I$ is a finitely presented $A'$-module, showing that $A'$ is coherent.
	\end{montobo}
	In \Cref{thm:4}, we give a sufficient list of axioms for a functor $\scr{F}$ to satisfy Gersten's injectivity for $\cal{S}_{\mathrm{val}}$. This list is entirely motivated by our proof, and in no way necessary. Demonstratively, even though the functors in \textsection\ref{section:toral_case} satisfy Gersten's injectivity, they violate our list of axioms. A functor $\scr{F}$ that satisfies our list of axioms shall be termed a $G$-theory-like functor (see \Cref{defn:G-theory_like_functor} below).
	\begin{montobo}[$G$-theory-like Functor]\label{montobo:K-theory-like-functor}
		Let $\rings$ be the category of commutative, unital rings, let $\scr{S}$ be a stable $\infty$-category and let $\scr{F}\colon\rings\to\scr{S}$ be a covariant functor. We suppose that $\scr{S}$ has a suitable notion of homotopy groups $\pi_q(-)$ which commute with filtered colimits, for all integer $q$. For example, $\scr{S}$ could be the $\infty$-category of `spectra'\footnote{Indeed, in this case, the homotopy groups are represented by suspensions of the `sphere spectrum' $\bb{S}$, which is a compact object (see \ccite{lurie-algebra}{Corollary~1.4.4.6}).} (see \ccite{lurie-algebra}{\textsection1.4.3}).
	\end{montobo}
		\bdf
			The functor $\scr{F}$ is said to satisfy \textit{localisation for a subcategory $\cal{S}\subseteq\rings$} if given any diagram \bud Z\arrow[r, hook, "i"]& X& U\arrow[hook']{l}[swap]{j}\eud of affine schemes such that $U=X\setminus Z$, $i$ is a finitely presented, closed immersion and $j$ is an open immersion between spectra of rings in $\cal{S}$, there exist
			\bun
				\item[$\circ$]  a pushforward morphism $i_{\ast}\colon\scr{F}(Z)\to\scr{F}(X)$, and
				\item[$\circ$]  an exact triangle \bud\scr{F}(Z)\arrow[r, "i_{\ast}"]&\scr{F}(X)\arrow[r,"j^{\ast}"]&\scr{F}(U).\eud
			\eun
		\edf
		A functor $\scr{F}$ that satisfies localisation for a category does not see nilpotents in the same. More precisely, such a functor $\scr{F}$ applied to the canonical morphism $A\twoheadrightarrow A_{\mathrm{red}}$, where the latter is the reduced ring of the former, is an isomorphism. We shall use this fact in the proof of \Cref{thm:4} without mention.  
		\par Given a small exact category $\cal{C}$, let $K(\cal{C})$ denote its `Quillen $K$-theory' spectrum\footnote{Although Quillen \ccite{quillen_k_theory}{\textsection2} initially defined the algebraic K-theory as a space using his $Q$-construction, Waldhausen \ccite{waldhausen_algebraic_k_theory_of_spaces}{1.5 and Appendix~1.9} demonstrated that the algebraic $K$-theory is actually an ``infinite loop space'' (in the sense of Adams \ccite{adams_infinite_loop_spaces}{\textsection1.4}). Consequently, according to op.~cit.~\textsection1.7, the algebraic $K$-theory can be regarded as a spectrum (cf.~\ccite{thomason_trobaugh}{1.5.2 and Theorem~1.11.2}).} defined in \ccite{quillen_k_theory}{\textsection2}. Given a ring $R$, let $K(R)$, termed its \textit{algebraic $K$-theory}, denote the Quillen $K$-theory of the category $\mathrm{Proj}(R)$ of finite type, projective $R$-modules. On the other hand, given a coherent ring $R$, let $G(R)$, termed its \textit{$G$-theory}, denote the Quillen $K$-theory of the category $\mathrm{Coh}(R)$ of $R$-finitely presented modules.
		\par By Quillen's d\'evissage theorem \ccite{quillen_k_theory}{\textsection5, Theorem~4}, it follows that the canonical inclusion $\mathrm{Proj}(R)\hookrightarrow\mathrm{Coh}(R)$ induces an isomorphism \bud\text{$K(R)\cong G(R)$ when $R\in\cal{S}_{\mathrm{val}}$},\eud i.e., when $R$ is a smooth algebra over a valuation ring. Indeed, for such a ring $R$, finitely presented $R$-modules admits a finite length resolution by finite type, projective $R$-modules. From the displayed isomorphism above, it follows that \bud\text{$G$-theory satisfies localisation for $\cal{S}_{\mathrm{val}}$ (see \ccite{antieau-mathew-morrow}{Proposition~2.5}).}\eud The rationale behind naming this section ``Gersten's Injectivity of $G$-theory'' instead of ``Gersten's Injectivity of $K$-theory'' becomes evident upon examining the proof of \Cref{thm:4}. It is apparent that $\mathrm{Coh}(R)$, rather than $\mathrm{Proj}(R)$, plays the key role in the proof.
		Let $\frak{Coh}$ be the $2$-category of additive categories of the form $\mathrm{Coh}(R)$, for some coherent ring $R$.
		\bdf
			A functor $\ca{K}\colon\frak{Coh}\to\scr{S}$ is said to \textit{satisfy the additivity property} if given any two coherent rings $R$ and $S$ as well as additive functors $f,g\colon\mathrm{Coh}(R)\to\mathrm{Coh}(S)$, we have an equivalence of functors $\ca{K}(f)+\ca{K}(g)\cong\ca{K}(f+g)$.
		\edf
		Quillen in \ccite{quillen_k_theory}{\textsection3, Corollary 1 to Theorem 2} shows that the algebraic $K$-theory satisfies the additive property.
		\bdf\label{defn:G-theory_like_functor}
			The functor $\scr{F}$ is called \textit{$G$-theory-like} if it satisfies the following properties, namely, if
			\bn[(i)] 
				\item\label{point:finite_presentation} it commutes with filtered colimits of rings,
				\item\label{point:localising_invariant} it satisfies localisation for $\cal{S}_{\mathrm{val}}$, and
				\item\label{point:Additivity_property} there exist a functor $\cal{K}$ that satisfies the additivity property and a functorial isomorphism $\scr{F}(R)\cong\mathcal{K}(\mathrm{Coh}(R))$, for any coherent ring $R$.
			\en
		\edf
		From the preceding remarks, it follows that $G$-theory is an example of a $G$-theory-like functor. For the subsequent discussion, we define the functors $\scr{F}_q(-)\colonequals\pi_q(\scr{F}(-))$, for each integer $q$.
		\begin{rem}~
			By our assumption, the homotopy functors $\pi_q(-)$ commute with filtered colimits, for all $q$. Consequently, if $\scr{F}(-)$ commutes with filtered colimits of rings, the functors $\scr{F}_q(-)$ will also commute with filtered colimits of rings, for all $q$.
		\end{rem}
	\par We are now ready to state and prove a generalisation of \Cref{thm:1.3}, which is recovered by plugging $\scr{F}(-)=K(-)$, the algebraic $K$-theory of rings.
	\bt\label{thm:4}
		Suppose that $\scr{F}$ is a $G$-theory-like functor.
		Given 
		\bun
			\item[$\circ$] a valuation ring $R$, 
			\item[$\circ$] an integral domain $\cal{A}$ that is a smooth $R$-algebra, 
			\item[$\circ$] a faithful $R$-algebra $A$ that is the semilocalisation of $\cal{A}$ at finitely many primes, and 
			\item[$\circ$] the subset $\ca{P}\subset\spec(A)$ of primes that correspond to generic points of the $R$-special fibre of $\spec(A)$,
		\eun the pullback morphism induces an injection $\scr{F}_i(A)\hookrightarrow \scr{F}_i(A_{\ca{P}})$, for all $i$.  
	\et
	\begin{rems}~
		\bn[(1)]
			\item Modulo the technicalities, the proof follows the arguments in \cite{gillet_levine}, which are, in itself, mixed-characteristic adaptations of Quillen's (\ccite{quillen_k_theory}{\textsection7, Theorem~5.11}). By an application of \Cref{defn:G-theory_like_functor}\eqref{point:finite_presentation}-\eqref{point:localising_invariant}, we reduce to show that $\scr{F}$ applied to $j\colon Y\hookrightarrow X$ is the zero morphism, where $X\colonequals\spec(\cal{A})$ and $Y$ is roughly the spreading out of $\spec(A/\frak{p})$. We then use Presentation Lemma~\ref{lem:presentation_lemma_Prufer} to further reduce to the case when $j$ is the vanishing locus of a principal ideal and has a retraction. Finally, in this case, we employ \Cref{defn:G-theory_like_functor}\eqref{point:Additivity_property} to conclude.
			\item We recall that the $R$-special fibre of $\spec(A)$ need not be connected, even when $R$ is a discrete valuation ring and $R\to\cal{A}$ is \'etale (see \Cref{rem:generic_fibre}\eqref{point:remark:generic_fibre}).
		\en
	\end{rems}
	\bs[Proof of \Cref{thm:4}]
		  Firstly, we shall reduce to the case when $R$ is of finite Krull dimension. Thanks to \Cref{lem:prufer_domain}\eqref{point:3_prufer_domain}, the ring $R$ can be written as filtered colimit of valuation rings $R_{\alpha}$ of finite Krull dimension. Since $\cal{A}$ is $V$-smooth, it is $V$-flat and $V$-finite type. Consequently, $\cal{A}$ is $R$-finitely presented (\ccite{raynaud_gruson}{Premi\`ere partie, théorème~3.4.6}). Therefore, possibly by reindexing, the smooth algebra $R\to\cal{A}$ descends to a smooth algebra $R_{\alpha}\to\cal{A}_{\alpha}$ (\stacks{0c0b}), for all $\alpha$. We suppose that $x_1,\ldots,x_n$ are the set of points corresponding to the semilocalisation $\mathcal{A}\to A$. Let $A_{\alpha}$ be the semilocalisation of $\cal{A}_{\alpha}$ at the images of the points $x_1,\ldots,x_n$ under the morphism $\spec(\cal{A})\to\spec(\cal{A}_{\alpha})$. Let $\frak{p}_{\alpha}\subset A_{\alpha}$ be a prime that corresponds to the generic point of the $R_{\alpha}$-special fibre of $\spec (A_{\alpha})$. Since the filtered colimit of $A_{\alpha}$ (respectively, $(A_{\alpha})_{\frak{p}_{\alpha}}$) is $A$ (respectively, $A_{\frak{p}}$) and $\scr{F}_i(-)$ commutes with filtered colimits of rings (\Cref{defn:G-theory_like_functor}\eqref{point:finite_presentation}), for all $i$, it suffices to show that the pullback morphism induces an injection $\scr{F}_i(A_{\alpha})\hookrightarrow \scr{F}_i((A_{\alpha})_{\frak{p}_{\alpha}})$, for each $\alpha$ and for all $i$. As a consequence, without loss of generality, we may assume that $R=R_{\alpha}$. In particular, \bud\text{we may assume that $R$ is of finite Krull dimension.}\eud Secondly, if $\cal{A}$ is of $R$-relative dimension $0$, then it is $R$-étale. In this case, $\cal{A}$ is a semilocal Pr\"ufer domain and $A=A_{\ca{P}}$. Thus, the claim is trivial. Therefore, we may assume that $\cal{A}$ is of $R$-relative dimension $d>0$.
		\par Thanks to \Cref{lem:snail_generic_point_of_special_fibre}\eqref{claim:2}, we have an equality \bud\textstyle A_{\ca{P}}=\colim_{f\in\mathrm{Flat}_{\ca{P}}(A/R)} A\inv{f}.\eud Similar to above, since $\scr{F}_i(-)$ commutes with colimits of rings, for any $i$, it suffices to show that $\scr{F}_i(A)\hookrightarrow \scr{F}_i(A\inv{f})$, for all $i$ and for each $f\in\mathrm{Flat}_{\ca{P}}(A/R)$. Since $\scr{F}$ satisfies localisation and it can be written as $\scr{F}=\scr{K}\circ\mathrm{Coh}$ (see \Cref{defn:G-theory_like_functor}\eqref{point:localising_invariant} and \eqref{point:Additivity_property}), it suffices to show that the morphism induced by applying $\scr{F}(-)$ to the pushforward morphism $\coh(A/(f))\to\coh(A)$ is zero. 
		Let $x_1,\ldots,x_n\in\spec(\cal{A})$ be the points corresponding to the semilocalisation $\mathcal{A}\to A$. Possibly by shrinking $\spec(\mathcal{A})$ to an affine neighbourhood of $x_1,\ldots,x_n$ and by spreading out, we may assume that $f\in\mathcal{A}$, whence the openness of the flat locus \stacks{0399} implies that, by further shrinking $\spec(\cal{A})$, we may suppose that $R\to\mathcal{A}/(f)$ is flat. Since, by assumption, $f$ does not vanish at the generic points of the $R$-special fibre of $\spec(\cal{A})$, the closed subscheme $i\colon Y\colonequals\spec(\mathcal{A}/(f))\hookrightarrow X\colonequals\spec\cal{A}$ does not contain any $R$-fibre of the target scheme (\stacks{0d4h}). Thus, applying Presentation Lemma~\ref{lem:presentation_lemma_Prufer}, at the cost of shrinking $X$ further, we get a smooth $R$-morphism $\pi\colon X\to\bb{A}^{d-1}_R$ of relative dimension $1$ such that $\pi|_Y$ is quasi-finite. Following the proof of \ccite{quillen_k_theory}{\textsection7, Theorem~5.11}, we define $C\colonequals Y\times_{\bb{A}^{d-1}_R}X$, i.e., we define $C$ via the pullback diagram given below. \begin{tpic}\label{diag:small_one}\node (C) at (0,0) {$C$}; \node (Y) at (0,-2) {$Y$}; \node (X) at (2,0) {$X$}; \node (A) at (2,-2) {$\bb{A}^{d-1}_R$}; \path[->] (X) edge node[right] {$\pi$} (A) (C) edge node[right] {$c$} (Y) edge node[above] {$q$} (X) (Y) edge node[below] {$\pi|_Y$} (A); \path[dashed, ->, bend left] (Y) edge node[left] {$s$} (C); \path[right hook->] (Y) edge node[below] {$i$} (X); \end{tpic} The section $s$, which is induced by $i$, is a closed immersion as $c$ is affine. However, by definition, $q$ is only quasi-finite, and not finite. This is the price we pay in mixed characteristic, as opposed to Quillen's presentation lemma \ccite{quillen_k_theory}{\textsection7, Lemma~5.12}, where $\pi|_Y$ can even be arranged to be finite. To deal with this shortcoming, for the rest of the proof we follow \ccite{gillet_levine}{\textsection2, proof of Theorem}. Thanks to \ccite{sga1}{corollaire~4.17}, the section $s$ of the smooth morphism $c$ is a regular immersion. Consequently, at the cost of shrinking $X$, we may assume that the ideal of vanishing of $s(Y)\subset C$ is principal, say with a generator $g$. By Zariski's Main Theorem \stacks{00qb}, there exists an open immersion $u\colon C\hookrightarrow\bar{C}$ with an extension of $q$, i.e., a finite morphism $\bar{q}\colon\bar{C}\to X$. We define $\bar{s}\colonequals u\circ s$. The diagram is depicted below.
		\begin{tpic}\label{diag:big_one}\node (C) at (0,0) {$C$}; \node (Y) at (-2.5,0) {$Y$}; \node (X) at (2.5,0) {$X$}; \node (Cbar) at (0,2) {$\bar{C}$}; \path[->] (C) edge node[above] {$q$} (X) (Cbar) edge node[above] {$\bar{q}$} (X); \path[bend right, ->] (C) edge node[above] {$c$} (Y); \path[right hook->] (Y) edge[bend right] node[above] {$i$} (X) edge node[above] {$\bar{s}$} (Cbar) edge node[below] {$s$} (C) (C) edge node[right] {$u$} (Cbar); \end{tpic}%
		We observe that $\bar{s}$ is a monomorphism that is also proper. Indeed, as $\bar{s}$ is a morphism from a proper $X$-scheme to a separated $X$-scheme, it is proper. Thus, it follows that $\bar{s}$ is a closed immersion (\stacks{04xv}). We need to show that $\scr{F}(-)$ applied to $i_{\ast}\colon\coh(Y)\to\coh(X)$ is the zero functor. Thanks to the fact that both $\bar{s}$ (since it's a closed immersion) and $\bar{q}$ are finite, it follows that each gives rise to a pushforward morphism between categories of coherent modules. Consequently, the equality $i=\bar{q}\circ\bar{s}$ induces an isomorphism of functors $i_{\ast}\cong\bar{q}_{\ast}\circ\bar{s}_{\ast}$. Therefore, it is enough to show that $\scr{F}(\bar{s}_{\ast})\cong 0$. On the other hand, the equality $c\circ s=\id_Y$ furnishes the isomorphism $s^{\ast}\circ c^{\ast}\cong\id_{\coh(Y)}$. As a consequence, it suffices to show that $\scr{F}(-)$ applied to $\bar{s}_{\ast}\circ s^{\ast}\colon\coh(C)\to\coh(\bar{C})$ is the zero functor. Let $Z\colonequals(\bar{C}\setminus C)$ be endowed with any closed subscheme structure. Admitting the fact that the ideal $\scr{I}$ of vanishing of $Z\subset\bar{C}$ is finitely generated (we prove this below in Claim 1), by localisation (\Cref{defn:G-theory_like_functor}\eqref{point:localising_invariant}), we have an exact triangle \bud\scr{F}(Z)\arrow[r]&\scr{F}(\bar{C})\arrow[r]&\scr{F}(C).\eud Therefore, it suffices to show that $\scr{F}(-)$ applied to $F\colonequals\bar{s}_{\ast}\circ s^{\ast}\circ u^{\ast}$ is the zero functor. We can further assume that the ideal $I$ of vanishing of $\bar{s}(Y)\subset\bar{C}$ is principal (we also prove this below in Claim 1), say with generator $\tilde{g}$. Since $\tilde{g}\in\scr{R}\colonequals\cal{O}_{\bar{C}}(\bar{C})$ is a nonzerodivisor, it induces an isomorphism $\scr{R}\cong I$ of $\scr{R}$-modules.
		\par Assuming $\bar{C}$ is coherent (we prove this in Claim 2), it follows that the exact sequence $0\to I\to\scr{R}\to\ca{A}/(f)\to 0$ induces a short exact sequence of functors $\coh(\bar{C})\to\coh(\bar{C})$ given by \bud0\tir&\id_{\coh(\bar{C})}\arrow[r,"\tilde{g}"]&\id_{\coh(\bar{C})}\tir&F\tir&0,\eud where we have invoked the isomorphism $\scr{R}\cong I$ to write the left morphism. Finally, thanks to the additivity (\Cref{defn:G-theory_like_functor}\eqref{point:Additivity_property}) of $\scr{F}(-)$, the above displayed exact sequence demonstrates that $\scr{F}(F)\cong 0$, as required. Thus, we are done modulo the claims made above, which we prove below. \vspace{0.15cm}
		\par Let us introduce some notations for convenience before proving the first claim. Specialising each $x_i$ to a closed point in $X$, we may, without loss of generality, assume that each $x_i$ is a closed point. Let $S\colonequals\{x_1,\ldots,x_n\}$, let $T\colonequals q^{-1}(S)$ and let $\bar{T}\colonequals\bar{q}^{-1}(S)$. Since $C\subset\bar{C}$, by definition, we have that $T\subset\bar{T}$. On the other hand, as $\bar{q}$ is finite, we conclude that $\bar{T}$ is a finite set.
		\vspace{0.2cm}
		\begin{claim}{~1}
			Without loss of generality, we may assume that both $\scr{I}$ and $I$ are principal.
		\end{claim}\vspace{0.1cm}
		\begin{claimproof}
			Firstly, we establish that $\scr{I}$ can be assumed to be principal. We semilocalise $\bar{C}$ at $\bar{T}$ (resp., $C$ at $T$) to obtain $\scr{R}_{\bar{T}}$ (resp., $\scr{R}_{T}$). Notably, $\spec(\scr{R}_{T})\subseteq\spec(\scr{R}_{\bar{T}})$ is an open subscheme. By prime avoidance \stacks{00ds}, we select an element $z\in\scr{I}\scr{R}_{\bar{T}}$ such that $z$ does not vanish at $ y_i$, for each $y_i\in T$. Consequently, according to the definition, we have \bud V(z)=Z\cap\spec(\scr{R}_{\bar{T}}).\eud Hence, by replacing $\scr{I}\scr{R}_{\bar{T}}$ with $z\scr{R}_{\bar{T}}$, we may assume that $\scr{I}\scr{R}_{\bar{T}}$ is principal. To show that $\scr{I}$ itself may be assumed to be principal, it is enough to spread out $z$ to a element of $\scr{R}$. By spreading out, there exist an open neighbourhood $U\subseteq\bar{C}$ of $\bar{T}$ and an element \bud\text{$\tilde{z}\in\mathcal{O}_{\bar{C}}(U)$ such that $V(\tilde{z})=Z\cap U$.}\eud We define $W\colonequals\overline{q}(\overline{C} \setminus U)$. We observe that $W \subseteq X$ is a closed subset since $\overline{q}$ is finite, and thus, a closed morphism. Then, we select a nonempty, affine open subset $X' \subseteq (X \setminus W)$. Finally, we replace diagram \eqref{diag:big_one} with its base change along $X'\hookrightarrow X$ (and by abuse of notation, refer to the objects with their respective names) to obtain an element $\tilde{z}\in\scr{R}$ whose vanishing locus is $Z$, as required. This completes the proof.
			\par It remains to show that $I$ can be assumed to be principal. Let $T_1\colonequals\bar{s}(Y)\cap\bar{T}$. In similar vein as the above, prime avoidance \stacks{00ds} furnishes an element $t\in I\scr{R}_{\bar{T}}$ such that $t$ does not vanish at $y_i$, for each $y_i\in(\bar{T}\setminus T_1)$. In consequence, by definition, we have \bud V(t)=Y\cap\spec(\scr{R}_{\bar{T}}).\eud Therefore, by substituting $I\scr{R}_{\bar{T}}$ with $t\scr{R}_{\bar{T}}$, we may assume that $I\scr{R}_{\bar{T}}$ is principal. A spreading out argument similar to above yields a nonempty, affine open subset $X'\subseteq X$. Upon base changing diagram \eqref{diag:big_one} along $X'\hookrightarrow X$ (and by abuse of notation, referring to the objects with their respective names), we obtain an element $\tilde{t}\in\scr{R}$ whose vanishing locus is $\bar{s}(Y)$, as required. Thus, the claim is proven.
			\par 
		\end{claimproof}\vspace{0.2cm}
		We prove the second claim below.
		\begin{claim}{~2}
			The scheme $\bar{C}$ is coherent.
		\end{claim} \vspace{0.1cm}
		\begin{claimproof}
			It is enough to show that $\bar{C}$ is $R$-finitely presented. By \ccite{raynaud_gruson}{Premi\`ere partie, théorème~3.4.6}, it suffices to verify that it is $R$-flat. Since $R$ is a valuation ring, it is equivalent to verify that $\bar{C}$ is $R$-torsion free. However, by the finer version of Zariski's Main Theorem \stacks{00qb}, we know that $\scr{R}\subseteq\scr{O}_C(C)$, where the larger ring is $R$-torsion free because $C$ is $R$-flat. Thus, the claim is proven.
		\end{claimproof}
			
	\es
	We are finally ready to establish the promised Gersten's injectivity in the case of $K$-theory (\Cref{thm:1.1}). By applying \Cref{thm:4} and employing standard reductions, the proof of the corollary below reduces to \ccite{kelly-morrow}{Theorem~3.1}.
	\bc\label{cor:Gersten_injectivity_for_valuation}
		Let $R$ be a semilocal Pr\"ufer domain and let $A$ be an integral domain that is the semilocalisation of a smooth $R$-algebra at finitely many primes. If $R$ contains a field, then the pullback morphism induces an injection $K_i(A)\hookrightarrow K_i(\Frac(A))$, for all $i$. 
	\ec
	\bs
		Similar to the beginning of the proof of \Cref{thm:4}, by a limit argument, we may, without loss of generality, assume that $R$ is of finite Krull dimension. Let $F\colonequals\Frac(R)$ and let $A_F\colonequals A\otimes_R F$. Thanks to \ccite{quillen_k_theory}{\textsection7, Theorem~5.11}, the restriction map induces an injection $K_i(A_F)\hookrightarrow K_i(\Frac(A))$, for each $i$. Consequently, it suffices to show that the restriction map induces an injection $K_i(A)\hookrightarrow K_i(A_F)$, for each $i$. 
		\par As a first case, we suppose that the claim is true when $R$ is a valuation ring. With this assumption, we shall show that the claim is true in general. To do so, we induct on \bud d(R)\colonequals\displaystyle\sum_{\substack{\fr{m}\in\mathrm{MaxSpec}(R)}}\dim(R_{\fr{m}}).\eud Let $\fr{m}\subset R$ be a maximal ideal and suppose that $a\in R$ is an element such that $V(a)=\{\fr{m}\}$ (such an $a$ exists thanks to prime avoidance \stacks{00ds}). Therefore, by definition, $\spec(R)\setminus\{\fr{m}\}=\spec(R\inv{a})$. Proceeding inductively, it suffices to show that $K_i(A)\hookrightarrow K_i(A\inv{a})$, for each $i$. Thanks to localisation \ccite{antieau-mathew-morrow}{Proposition~2.5}, we obtain a morphism of exact triangles \bud G(A/(aA))\arrow[d]\tir& K(A)\tir\arrow[d]& K(A\inv{a})\arrow[d]\\ G(A_{\fr{m}}/(aA_{\fr{m}}))\tir& K(A_{\fr{m}})\tir& K(A_{\fr{m}}\inv{a}).\eud However, the canonical morphism induces an isomorphism $A/(aA)\cong A_{\fr{m}}/(aA_{\fr{m}})$, demonstrating that the right square is Cartesian. As $R_{\fr{m}}$ is a valuation ring, by our assumption, the lower left horizontal morphism is zero. Consequently, it follows that the upper left horizontal morphism is also zero, establishing the required injection. Thus, in this case, the proof is complete.
		\par As a consequence, it remains to establish the claim when $R$ is a valuation ring. Let $\ca{P}\subset\spec(A)$ be the subset of primes that correspond to the generic points of the $R$-special fibre of $\spec(A)$. By \Cref{thm:4}, we have $K_i(A)\hookrightarrow K_i(A_{\cal{P}})$, for each $i$. Since $\Frac(A_{\ca{P}})=\Frac(A)$, without loss of generality, we can assume that $A=A_{\ca{P}}$. Moreover, under this assumption, thanks to \Cref{lem:snail_generic_point_of_special_fibre}\eqref{claim:1}, the ring $A$ is a Pr\"ufer domain. Using an induction argument similar to the one above, we reduce the problem to showing the claim when $A$ is a valuation ring. In this case, the required injectivity is the content of \ccite{kelly-morrow}{Theorem~3.1}.
	\es
	\br
		In effect, the proof of \Cref{cor:Gersten_injectivity_for_valuation} demonstrates that if Gersten's injectivity holds for all valuation rings, it extends to all smooth algebras over valuation rings. However, the mixed-characteristic counterpart of \ccite{kelly-morrow}{Theorem~3.1} remains largely unexplored. Nevertheless, there are some preliminary results in the case of perfectoid valuation rings (see \ccite{antieau-mathew-morrow}{Theorem~1.2}).
	\er
	\section{Toral Case of Gersten's Injectivity}\label{section:toral_case}
	In this section, unless otherwise stated, all the cohomology groups that appear are \'etale. Let $X$ be a `non-singular' scheme and let $j\colon U\hookrightarrow X$ be an open subscheme whose complement is of depth $\ge 2$. Consider an isotrivial $X$-torus $T$. The primary objective of this section is to establish \Cref{lem:Brauer_group_injects}. Our key technical tool is a version of purity of torsors under tori (see \eqref{diag:lifting_torus_torsor_S2}), which essentially provides a necessary and sufficient condition for a $T$-torsor on $U$ to extend globally over $X$. We prove this purity via an application of a weak version of Auslander--Buchsbaum formula \eqref{pushforward:picard_group}, which demonstrates that $j_{\ast}$ preserves reflexive sheaves.

	\par To start, we recall the notion of reflexive sheaves.
	\begin{montobo}[Reflexive sheaves] 
		Let $X$ be a scheme. The \textit{dual} $\scr{F}^\vee$ of an $\ca{O}_X$-module $\scr{F}$ is defined to be the $\ca{O}_X$-module $\underline{\Hom}_{\ca{O}_X}(\scr{F},\ca{O}_X)$. A coherent $\ca{O}_X$-module $\scr{F}$ is called \textit{reflexive} if for every $x\in X$, there is a neighbourhood $x\in U\subset X$ such that the canonical morphism \bud\beta_{\scr{F}|_U}\colon\scr{F}|_U\isom&\scr{F}|_U^{\vee\vee}\text{\hspace{0.5cm}is an isomorphism.}\eud Given a coherent $\ca{O}_X$-module $\scr{F}$ and a presentation \bd\label{presentation:reflexive_sheaf}\ca{O}_X^{\oplus m}\tir&\ca{O}_X^{\oplus n}\tir&\scr{F}\tir&0,\ed we can dualise to obtain a short exact sequence \bd\label{dual:presentation_reflexive_sheaf}0\tir&\scr{F}^{\vee}\tir&\ca{O}_X^{\oplus n}\tir&\ca{O}_X^{\oplus m}.\ed Therefore, by \cite[\href{https://stacks.math.columbia.edu/tag/01BY}{Tag 01BY}]{stacks-project}, if $X$ is a locally coherent scheme, then for a coherent sheaf (in particular, reflexive) $\scr{F}$, the dual $\scr{F}^{\vee}$ is coherent.
	\end{montobo} 
	
	In preparation for \Cref{prop:extend_tori_torsor}, we prove the following lemma, drawing inspiration from \ccite{gabber-ramero-foundations}{Proposition~11.3.8} and \ccite{ct-sansuc_extend_torosor_tori}{Lemma 2.1}.
	\bl\label{lem:Gabber_Ramero_S2_condition}
	For a locally coherent scheme $X$, a quasi-compact open $j\colon U\hookrightarrow X$ such that at each point $z\in Z\colonequals X\setminus U$, we have\footnote{A module $M$ over a local ring $(A,\fr{m})$ has $\dep_A(M)\ge d$, if there is an $M$-regular sequence $x_1,\ldots,x_d\in\fr{m}$; the depth of $A$ is $\dep_A(A)$ (see \ccite{egaIV_1}{Chapitre 0, Définition~15.1.7 and \textsection15.2.2)}). There is no condition on the quotients being nonzero.} $\dep(\ca{O}_{X,z})\ge 2$,
	and a reflexive $\ca{O}_X$-module $\scr{F}$, the restriction induces an isomorphism \bd\label{diag:lifting_line_bundles_S2}\scr{F}\isom& j_{\ast}j^{\ast}\scr{F}.\ed 
	Moreover, if $X$ is reduced, for a reflexive $\ca{O}_U$-module $\scr{G}$, \bd\label{pushforward:reflexive_module}\text{the pushforward $j_{\ast}\scr{G}$ is a reflexive $\ca{O}_X$-module.}\ed
	\el
	\bs
	\eqref{diag:lifting_line_bundles_S2}$\colon$ Thanks to \ccite{ces-scholze}{Lemma 7.2.7(b)}, the restriction induces an isomorphism \bd\label{diag:lifting_functions_S2}\ca{O}_X\isom& j_{\ast}\ca{O}_U.\ed We shall reduce to the special case $\scr{F}=\ca{O}_X$. Since it is enough to show \eqref{diag:lifting_line_bundles_S2} locally, given a reflexive $\ca{O}_X$-module $\scr{F}$, we may assume that there is a presentation \eqref{presentation:reflexive_sheaf} of $\scr{F}^{\vee}$, which can be dualised to obtain a short exact sequence like \eqref{dual:presentation_reflexive_sheaf}. Since $j_\ast$ is left exact, this gives us a commutative diagram \bud0\tir&\scr{F}\arrow[d]\tir&\ca{O}_X^{\oplus n}\arrow{d}{\sim}[swap]{\eqref{diag:lifting_functions_S2}}\tir&\ca{O}_X^{\oplus m}\arrow{d}{\eqref{diag:lifting_functions_S2}}[swap]{\sim}\\0\tir&j_{\ast}j^{\ast}\scr{F}\tir&\ca{O}_X^{\oplus n}\tir&\ca{O}_X^{\oplus m}, \eud from which we are reduced to the case when $\scr{F}=\ca{O}_X$, and we are done.\vspace{0.25cm}
	\par\eqref{pushforward:reflexive_module}$\colon$ In view of \eqref{diag:lifting_line_bundles_S2}, it is enough to show that there exists a reflexive $\ca{O}_X$-module $\scr{F}$ such that $\scr{F}\mid_U=\scr{G}$. By \cite[\href{https://stacks.math.columbia.edu/tag/0G41}{Tag 0G41}]{stacks-project}, there is a finitely presented $\ca{O}_X$-module $\scr{F}'$ such that $\scr{F}'\mid_U=\scr{G}$. Thanks to \cite[\href{https://stacks.math.columbia.edu/tag/01BZ}{Tag 01BZ}]{stacks-project}, keeping in mind that $X$ is locally coherent, the $\ca{O}_X$-module $\scr{F}'$ is automatically coherent. Taking $\scr{F}=\scr{F}'^{\vee\vee}$, we note that $\scr{F}\mid_U=\scr{F}'^{\vee\vee}\mid_U=\scr{G}^{\vee\vee}=\scr{G}$ (using the fact that $\scr{G}$ is reflexive). It remains to check that $\scr{F}$ is reflexive, for which we follow the proof of \cite[\href{https://stacks.math.columbia.edu/tag/0AY4}{Tag 0AY4}]{stacks-project}. Since the result is local, it can be assumed that $X=\spec A$ is affine. Choosing a presentation \bud A^{\oplus m}\to A^{\oplus n}\to\Gamma(\spec A,\scr{F'})\to 0,\eud and dualising it, in order to conclude, it is sufficient to show the following claim.\vspace{0.1cm}
	\begin{claim}{}
		Given an exact sequence \bud0\to M\to M'\to M''\eud of finitely presented $A$-modules, the module $M$ is reflexive if $M'$ and $M''$ are reflexive.
	\end{claim}\vspace{0.1cm}
	\begin{claimproof}
		We suppose that $M'$ and $M''$ are reflexive. Proceeding as in the proof of \cite[\href{https://stacks.math.columbia.edu/tag/0EB8}{Tag 0EB8}]{stacks-project}, we shall show that $M$ is reflexive. Double dualising the displayed short exact sequence in the claim and writing down canonical morphisms, we get the following morphism of complexes \bud M\tir\arrow[d]&M'\tir\arrow[d]&M''\arrow[d]\\\Hom_A(\Hom_A(M,A),A)\arrow[r,"\alpha"]&\Hom_A(\Hom_A(M',A),A)\tir&\Hom_A(\Hom_A(M'',A),A).\eud By the assumption, the middle and the right vertical arrows are isomorphisms. We need to show that the left vertical arrow is an isomorphism. It suffices to show that $\alpha$ is injective. We consider module $Q$ defined by the exact sequence $\Hom_A(M',A)\to\Hom_A(M,A)\to Q\to 0$. Letting $K$ be the total ring of fractions of $A$ (see \cite[\href{https://stacks.math.columbia.edu/tag/00EW}{Tag 00EW}]{stacks-project}), by the finite presentation property \cite[\href{https://stacks.math.columbia.edu/tag/0583}{Tag 0583}]{stacks-project}, tensoring the exact sequence with $K$, we obtain the exact sequence \bud\Hom_K(M'\otimes_A K,K)\to\Hom_K(M\otimes_A K,K)\to Q\otimes_A K\to 0.\eud However, since $K$ is a product of fields, the injection $M\otimes_A K\hookrightarrow M'\otimes_A K$ is split, consequently, $Q\otimes_A K=0$, implying that $Q$ is a torsion $A$-module. In that case, $\Hom_A(Q,A)=0$, ensuring that $\alpha$ is injective.
	\end{claimproof}
	\es


	Given a flat (resp., smooth) group scheme $G$ over a scheme $S$ and an $S$-scheme $S'$, we let $\mathbf{B}G(S')$ denote the category of fppf locally (resp., \'etale locally) trivial $G$-torsors on $S'$. Likewise, let $\mathbf{B}G$ be the presheaf on the category of $S$-schemes defined by $S'\mapsto\mathbf{B}G(S')$ (see \cite[\href{https://stacks.math.columbia.edu/tag/0048}{Tag 0048}]{stacks-project}).
	\par We are now ready to prove our version of the Auslander--Buchsbaum formula \eqref{pushforward:picard_group}. As a consequence, we derive a categorical variant of the purity of torsors under tori \eqref{diag:lifting_torus_torsor_S2}.
	\bp\label{prop:extend_tori_torsor} 
	Let $R$ be a Pr\"ufer domain, let $X$ be a smooth, integral $R$-scheme and let $j\colon U\hookrightarrow X$ be a quasi-compact open such that at each point $x\in Z\colonequals X\setminus U$ with $f(x)=y$, we have $\dim(\ca{O}_{X_{y},x})+\min(1,\dim(R_y))\ge 2$.
	Then, for a locally free $\ca{O}_U$-module $\scr{L}$ of rank $1$,\bd\label{pushforward:picard_group}\text{the pushforward $j_{\ast}\scr{L}$ is a locally free $\ca{O}_X$-module of rank $1$},\ed
	in particular, for any \'etale $X$-scheme $X'$, the restriction induces an equivalence of categories
	\bd\label{diag:lifting_picard_groups_S2}\mathbf{B}\bb{G}_m(X')\isom&\mathbf{B}\bb{G}_m(X'\times_X U).\ed
	More generally, for an $X$-torus $T$ and for any \'etale $X$-scheme $X'$, the restriction induces an equivalence of categories \bd\label{diag:lifting_torus_torsor_S2}\mathbf{B}T(X')\isom&\mathbf{B}T(X'\times_X U),\ed 
	in particular, \bd\label{diag:downgraded_version_lifting_torus_torsor_S2} H^q(X', T)\cong H^q(X'\times_X U, T),\text{\hspace{0.25cm}for }q\le 1.\ed
	\ep
	\br
		The seemingly unusual inequality $\dim(\ca{O}_{X_{y},x})+\min(1,\dim(R_y))\ge 2$ serves a specific purpose: ensuring that $\depth(\ca{O}_{X,x})\ge 2$ at each point $x\in Z$. Relatedly, Pr\"ufer domains are far from being Cohen--Macaulay. In effect, as previously noted in \textsection\ref{section:Prufer}, by prime avoidance \stacks{00ds}, for any valuation ring $V$ of finite Krull dimension, there exists an element $a\in V$ so that $V(a)=\{\fr{m}\}$, where $\fr{m}\subset V$ is the maximal ideal.  
	\er
	\bs[Proof of \Cref{prop:extend_tori_torsor}]
	We show that $\dep(\ca{O}_{X,z})\ge 2$, at any point $z\in Z$. Let $f\colon X\to\spec R$. Thanks to \ccite{egaIV_3}{Th\'eor\`eme 11.3.8} (specifically, (c)$\implies$(a)), it suffices to argue that $\dep(\ca{O}_{X_{f(z)},z})+\min(1,\dim(R_{f(z)}))\ge 2$. It follows from the hypothesis and the equality $\dep(\ca{O}_{X_{f(z)},z})=\dim(\ca{O}_{X_{f(z)},z})$, which is true because $f$ is smooth.
	\par We show the key claim, i.e., \eqref{pushforward:picard_group}, below. The claim \eqref{diag:lifting_picard_groups_S2} is a consequence of \eqref{diag:lifting_line_bundles_S2} and \eqref{pushforward:picard_group}. Since the torus $T$ trivialises \'etale locally on $X$, the claims \eqref{diag:lifting_torus_torsor_S2} and \eqref{diag:downgraded_version_lifting_torus_torsor_S2} reduce to \eqref{diag:lifting_picard_groups_S2} by an \'etale descent argument.
	\par\eqref{pushforward:picard_group}$\colon$ We follow the proof of \ccite{gabber-ramero-foundations}{Proposition 11.4.1(iv)}. A complex of sheaves on $X$ of abelian groups that is concentrated in cohomological degree $0$ with the $0$-th term $\scr{A}$ shall be denoted by $\scr{A}[0]$. Thanks to \eqref{pushforward:reflexive_module}, the pushforward $\scr{M}\colonequals j_{\ast}\scr{L}$ is a reflexive $\ca{O}_X$-module, and also, torsion-free and hence, flat over $R$. Assuming that $\scr{M}[0]$ is a perfect $\ca{O}_X$-complex, the result follows. Indeed, letting $\det(\scr{M}[0])$ be the determinant line bundle of $\scr{M}[0]$ (see \ccite{knudsen_mumford}{Theorem 1}), there is a sequence of isomorphisms \bud \scr{M}\isom&j_{\ast}\scr{L}\isom&j_{\ast}\det (\scr{L}[0])\isom&j_{\ast}j^{\ast}(\det(\scr{M}[0]))&\det(\scr{M}[0])\arrow{l}{\sim}[swap]{\eqref{diag:lifting_line_bundles_S2}},\eud from which the result follows. We verify that $\scr{M}[0]$ is a perfect $\ca{O}_X$-complex. It suffices to assume that $X=\spec A$ is affine and to show that $\scr{M}[0]$ is quasi-isomorphic to a bounded complex of finite free $A$-modules (\cite[\href{https://stacks.math.columbia.edu/tag/0BCJ}{Tag 0BCJ}]{stacks-project}). Let $M\colonequals\Gamma(A,\scr{M})$. By \cite[\href{https://stacks.math.columbia.edu/tag/0G9A}{Tag 0G9A}]{stacks-project}, it suffices to show that \bd\label{displayed_above}\text{$\Ext^q_A(M,N)=0$ for any finitely presented $A$-module $N$ and any $q\gg 0$}\ed (the complex $M[0]$ is pseudo-coherent because $M$ is a coherent $A$-module). Using the finitely presented property of the variable $N$ in \eqref{displayed_above}, it is enough to show that there exists an integer $n$ such that $\mathrm{proj.dim}_{A_{\fr{p}}}(M_{\fr{p}})\le n$ for any prime $\fr{p}\subset A$ (see \ccite{weibel}{Lemma~3.3.8}). This follows from \ccite{gabber-ramero-foundations}{Proposition 11.4.1(ii)} (or \ccite{ning_valuation}{Lemma~7.2(i)}), which shows that taking $n=\dim_R(A)$ suffices. Hence, we are done.

	\vspace{0.25cm}\par\eqref{diag:lifting_picard_groups_S2}$\colon$ \comment{If $x$ is not contained in the image of $X'\to X$, the statement is trivial. Henceforth, we assume that $x$ lies in the image of $X'\to X$. }Let $j'\colon U'\colonequals X'\times_X U\hookrightarrow X$ be the inclusion. Since $\dim_X(X')=0$, we obtain that $(X',U')$ satisfies the hypothesis of \Cref{prop:extend_tori_torsor}. We shall show that the pushforward $j'_{\ast}\colon\mathbf{B}\bb{G}_m(U')\to\mathbf{B}\bb{G}_m(X')$, which is well defined as a consequence of \eqref{pushforward:picard_group}, is an inverse to the pullback $j'^{\ast}\colon\mathbf{B}\bb{G}_m(X')\to\mathbf{B}\bb{G}_m(U')$. Since the equality $j'^{\ast}j'_{\ast}=\id$ follows from the definition, it suffices to show that for a line bundle $\scr{L}\in\mathbf{B}\bb{G}_m(X')$, the restriction induces an isomorphism $\scr{L}\iso j'_{\ast}j'^{\ast}\scr{L}$. However, this results from \eqref{diag:lifting_line_bundles_S2}.\vspace{0.25cm} 
	\par\eqref{diag:lifting_torus_torsor_S2}$\colon$ Let $U'\colonequals X'\times_X U$. Thanks to \cite[\href{https://stacks.math.columbia.edu/tag/04UK}{Tag 04UK}]{stacks-project}, $\mathbf{B}T_{X'}$ satisfies \'etale descent, consequently, the same holds for the presheaf $j_{\ast}\mathbf{B}T_{U'}$. By \ccite{sga3ii}{Exposé~X, Corollaire~4.5}, there exists an \'etale surjection $\tilde{X}\to X'$ that splits $T$. In view of the \'etale descent property of $\mathbf{B}T_{X'}$ and $j_{\ast}\mathbf{B}T_{U'}$, it suffices to show that for any \'etale $\tilde{X}$-scheme $X''$, the restriction induces an equivalence of categories $\mathbf{B}T_{\tilde{X}}(X'')\iso \mathbf{B}T_{\tilde{X}}(X''\times_{X'} U')$, which is the content of \eqref{diag:lifting_picard_groups_S2}.\vspace{0.25cm}
	\par\eqref{diag:downgraded_version_lifting_torus_torsor_S2} This follows from \eqref{diag:lifting_torus_torsor_S2}. Indeed, $H^0(X',T)$ (resp., $H^0(U', T)$) is the automorphism group of the trivial torsor in $\mathbf{B}T(X')$ (resp., in $\mathbf{B}T(U')$) and $H^1(X',T)$ (resp., $H^1(U',T)$) is the group of isomorphism classes of objects in $\mathbf{B}T(X')$ (resp., in $\mathbf{B}T(U')$).
	\es
	With \Cref{prop:extend_tori_torsor} proven, to state our main theorem, we now need to recall the notion of flasque tori.
	\begin{montobo}[Flasque Torus and Flasque Resolution]\label{montobo:flasque_torus_resolution}
		We recall some definitions from \ccite{ct-sansuc_flasque_tori}{\textsection0.5}. Let $G$ be a finite group. A finitely generated, free $\bb{Z}$-module $\ca{P}$ with a linear action of $G$ is called a \textit{permutation module} if $\ca{P}$ admits a $G$-stable $\bb{Z}$-basis. A finitely generated, free $\bb{Z}$-module $\ca{F}$ with a linear action of $G$ is called \textit{flasque} if \bud \Ext^1_{\bb{Z}[G]}(\ca{F},\ca{P})=0 \text{\hspace{0.5cm}or, equivalently,\hspace{0.5cm}} H^1(G,\Hom_{\bb{Z}}(\ca{F},\ca{P}))=0,\eud for any permutation $\bb{Z}[G]$-module $\ca{P}$. For example, the trivial $\bb{Z}[G]$-module $\bb{Z}$ is a permutation module, and as a consequence, for a flasque $\bb{Z}[G]$-module $\ca{F}$, we have $\Ext^1_{\bb{Z}[G]}(\ca{F},\bb{Z})=0$. 
		\par Let $X$ be a scheme. An $X$-torus $T$ is called \textit{isotrivial} if it is split by a finite \'etale surjection $\tilde{X}\to X$. The \textit{character group} of an $X$-torus $T$ is the sheaf of abelian groups $T^{\vee}\colonequals\underline{\Hom}_{X\text{-gps}}(T,\bb{G}_{m,X})$. An isotrivial $X$-torus $T$ is called \textit{quasi-trivial} (resp., \textit{flasque}) if for any connected component $Z\subset X$, there exists a connected, Galois, finite \'etale cover $\tilde{Z}\to Z$ that splits $T$ such that the induced $\bb{Z}[\Gal(\tilde{Z}/Z)]$-module $T^{\vee}(\tilde{Z})$ is a permutation module (resp., is flasque) (see \ccite{ct-sansuc_flasque_tori}{Definition 1.2}). In fact, by op.~cit.~Lemma 1.1, any connected, Galois, finite \'etale cover $\tilde{Z}\to Z$ that splits $T$ can be chosen in the previous definition, and the choice of the closed subscheme structure on $Z\subset X$ is irrelevant. As we might expect, the notions of flasque and quasi-trivial tori are preserved under base change (cf., op.~cit.~Proposition 1.3). For a connected scheme $X$, any quasi-trivial torus $Q$ can be written as a finite product of Weil restrictions $\mathrm{Res}_{X_i/X}(\bb{G}_{m,X_i})$, for finite \'etale covers $X_i\to X$ (see \ccite{ces_torsors}{Lemma~A.2.6}).
		\par Thanks to \ccite{ct-sansuc_flasque_tori}{Proposition 1.3}, given an isotrivial torus $T$ on a scheme $X$ whose connected components are open (for example, a scheme which has finitely many connected components, like the spectrum of a semilocal ring), there exists two different \textit{flasque resolutions} of $T$, namely exact sequences \bd\label{exact_sequence:flasque-resolution} 1\to F\to Q\to T\to 1,\text{ where $F$ is a flasque and $Q$ is a quasi-trivial $X$-torus, and}\ed \bd\label{exact_sequence:flasque-resolution2} 1\to T\to F'\to Q'\to 1,\text{ where $F'$ is a flasque and $Q'$ is a quasi-trivial $X$-torus.}\ed
	\end{montobo}
	\vspace{0.25cm}
	We are now prepared to prove our main theorem. Before stating the theorem, let us outline the steps of its proof. To show Gersten's injectivity for a flasque $A$-torus $F$, standard reductions show that it is enough to establish vanishing \bd\label{vanishing_intro}H^2_{\{z\}}(A_z, F)=0\ed of the local cohomology at a point $z\in\spec(A)$. If $z$ is of codimension $1$ or if it is the generic point of an $R$-fibre of $\spec(A)$, \Cref{lem:snail_generic_point_of_special_fibre} shows that $A_{z}$ is a valuation ring. Therefore, in this case, \eqref{vanishing_intro} follows from the results in \ccite{ning_valuation}{\textsection2}. Otherwise, our analysis shows that it follows from \Cref{prop:extend_tori_torsor}.
	\bt\label{lem:Brauer_group_injects}
	Given a semilocal Pr\"ufer domain $R$, an integral domain $A$ that is $R$-essentially smooth, a quasi-compact open $U\hookrightarrow\spec A$, and a flasque $A$-torus $F$, \vspace{0.2cm}
	\bn[(i)]
		\item\label{theta1F_surjection} the morphism $H^1(A,F)\to H^1(U,F)$ is surjective, and\vspace{0.2cm}
		\item\label{theta2Finjection} the morphism $H^2(A,F)\to H^2(U,F)$ is injective.
	\en
	\et
	\bs
	We reduce to showing \eqref{theta1F_surjection} and \eqref{theta2Finjection} for integral domains that are $R$-smooth. Let $\ca{A}$ be an $R$-smooth, integral domain such that $A$ is a semilocalisation of $\ca{A}$. 
	Assuming that \Cref{lem:Brauer_group_injects} holds for integral domains of the form $\ca{A}\inv{f}$, for some $f\in\ca{A}$, a limit argument and the facts that \'etale cohomology commutes with filtered colimits of rings (see \cite[\href{https://stacks.math.columbia.edu/tag/09YQ}{Tag 09YQ}]{stacks-project}) and that colimits commute with cokernels will then show that \eqref{theta1F_surjection} is true for $A$. In a similar vein, we reduce to showing \eqref{theta2Finjection} for rings of the form $\ca{A}\inv{f}$, for some $f\in\ca{A}$. Thus, without loss of generality, we assume that $A$ is $R$-smooth.
	\par Again, thanks to \Cref{lem:prufer_domain}\eqref{point:3_prufer_domain}, by a limit argument (see \cite[\href{https://stacks.math.columbia.edu/tag/09YQ}{Tag 09YQ}]{stacks-project}), we may assume that $R$ has finite Krull dimension.
		Letting  $Z\colonequals\spec A\setminus U$, an analysis of the long exact sequence of cohomology with supports
		\begin{tpic*}\label{long_exact_sequence:cohomology_with_supports} \node (d1) at (0,0) {$\cdots$}; \node (h1r) at (2,0) {$H^1(A,F)$}; \node (h1u) at (6,0) {$H^1(U,F)$}; \node (h2mr) at (-0.75,-1.25) {$H^2_{Z}(A,F)$}; \node (h2r) at (2,-1.25) {$H^2(A,F)$}; \node (h2u) at (6,-1.25) {$H^2(U,F)$}; \node (d2) at (9,-1.25) {$\cdots$}; \path[->,font=\scriptsize,>=angle 90] (d1) edge (h1r) (h1r) edge (h1u) (h2mr) edge (h2r) (h2r) edge (h2u) (h2u) edge (d2); \draw[->,gray]
			(h1u) edge[out=-9,in=170] (h2mr);
		\end{tpic*}indicates that it suffices to establish that $H^2_Z(A,F)=0$. On the other hand, since the $R$-fibres of $\spec(A)$ are Noetherian, the topological space $\spec(A)$ itself is Noetherian. Therefore, thanks to the coniveau spectral sequence associated to the filtration by supports inside closed subschemes \ccite{ilo14}{Expos\'e XVIII-A, \textsection2.2.1} (see also \ccite{brauerIII}{Section 10.1}) \bd\label{sseq:coniveau_1} E_1^{p,q}\colon\bigoplus_{z\in Z\text{ with } \dim A_{z}=p}H^{p+q}_{\{z\}}(A_z,F)\Rightarrow H^{p+q}_{Z}(A,F),\ed we are reduced to establishing that \bd\label{vanishing}\text{$H^{q}_{\{z\}}(A_{z},F)=0$, for each $z\in Z$ and for all $q\le 2$.}\ed We first examine the following two cases, which readily follow from results in the literature.
		\bun 
			\item[$\circ$] If $z$ is a generic point of an $R$-fibre, then \Cref{lem:snail_generic_point_of_special_fibre} shows that $A_{z}$ is a valuation ring. Indeed, letting $x$ be the image of $z$ along $\spec A\to\spec R$, the image of $z$ in $\spec A_x$ is a generic point of the special fibre over the valuation ring $R_x$. In this case, to prove \eqref{vanishing}, we again write a long exact sequence of cohomology with supports. Ultimately, it follows from \ccite{ning_valuation}{Lemma 2.3, Proposition~2.4 and Corollary~2.5}. 
			\item[$\circ$] If $z$ is a point of codimension $1$ in the generic fibre, then $A_z$ is an equicharacteristic discrete valuation ring. In a similar vein as above, vanishing \eqref{vanishing} follows.
		\eun
		Consequently, without loss of generality, we may assume that $z$ is neither a generic point of an $R$-fibre of $\spec(A)$ nor it is a point of codimension $1$ in the $R$-generic fibre. This assumption is important for our eventual application of \Cref{prop:extend_tori_torsor}. Thanks to the local-to-global spectral sequence (\ccite{sga4ii}{Expos\'e V, Proposition 6.5}) \bud E_2^{p,q}\colon H^p(A_z,\scr{H}^q_{\{z\}}(A_z,F))\Rightarrow H^{p+q}_{\{z\}}(A_z,F),\eud it is enough to demonstrate that $\scr{H}^{q}_{\{z\}}(A_z,F)=0$, for all $q\le 2$. Letting $A_{\overline{z}}$ be the strict Henselisation of $A_{z}$ and by taking stalks, it is equivalent to show that $H^{q}_{\{\overline{z}\}}(A_{\overline{z}}, F)=0$, for all $q\le 2$. An inspection of the long exact sequence of cohomology with supports \begin{tpic*}\label{long_exact_sequence:cohomology_with_supports_local} \node (d0) at (-2,0) {0}; \node (d1) at (0,0) {$H^0_{\{\overline{z}\}}(A_{\overline{z}},F)$}; \node (h1r) at (3,0) {$H^0(A_{\overline{z}},F)$}; \node (h1u) at (7,0) {$H^0(\spec A_{\overline{z}}\setminus\{\overline{z}\},F)$}; \node (h2mr) at (0,-1.375) {$H^1_{\{\overline{z}\}}(A_{\overline{z}},F)$}; \node (h2r) at (3,-1.375) {$H^1(A_{\overline{z}},F)$}; \node (h2u) at (7,-1.375) {$H^1(\spec A_{\overline{z}}\setminus\{\overline{z}\},F)$}; \node (h3r) at (0,-2.75) {$H^2_{\{\overline{z}\}}(A_{\overline{z}},F)$}; \node (h3u) at (3,-2.75) {$H^2(A_{\overline{z}},F)$}; \node (d3) at (6,-2.75) {$\cdots$}; \path[->,font=\scriptsize,>=angle 90] (d0) edge (d1) (d1) edge (h1r) (h1r) edge (h1u) (h2mr) edge (h2r) (h2r) edge (h2u) (h3r) edge (h3u) (h3u) edge (d3); \draw[->,gray]
			(h1u) edge[out=-9,in=170] (h2mr) (h2u) edge[out=-9,in=170] (h3r);
		\end{tpic*}and the vanishing of the positive degrees of \'etale cohomology of strictly Henselian rings \cite[\href{https://stacks.math.columbia.edu/tag/03QO}{Tag 03QO}]{stacks-project} proves that \bd\label{diag:local_vanishing}H^2_{\{\overline{z}\}}(A_{\overline{z}},F)\cong H^1(\spec A_{\overline{z}}\setminus\{\overline{z}\},F).\ed In conclusion, it is left to show that $H^{q}(\spec A_{\overline{z}},F)\cong H^q(\spec A_{\overline{z}}\setminus\{\overline{z}\},F)$, for $q=0$ and $q=1$. Since \'etale cohomology commutes with cofiltered limits of schemes, by the definition of strict Henselisation, it remains to apply \Cref{prop:extend_tori_torsor} to $X=\spec A, Z=\overline{\{z\}}$ and $T=F$ (since the topological space $\spec A$ is Noetherian, any open subset of $\spec A$ is quasi-compact).\es

				Finally, we are ready to establish Gersten's injectivity in the case of $H^1(-, T)$ as well as $H^2(-,T)$, for a torus $T$ over a smooth algebra over a semilocal Pr\"ufer domain. In other words, we prove a generalisation of \Cref{thm:1.2} below. By employing flasque resolutions of isotrivial tori, we derive this result as a consequence of \Cref{lem:Brauer_group_injects}.
				\bc\label{cor:torus_Gersten_injectivity}
				For a semilocal Pr\"ufer domain $R$, a semilocal, integral domain $A$ that is $R$-essentially smooth and an $A$-torus $T$,\vspace{0.2cm}
				\bn[(i)]
					\item\label{theta1Tinjection} the morphism $H^1(A,T)\hookrightarrow H^1(K,T)$ is injective, and\vspace{0.2cm}
					\item\label{theta2Tinjection} the morphism $H^2(A,T)\hookrightarrow H^2(K,T)$ is injective.
				\en
				\ec
				\bs
				\eqref{theta1Tinjection}$\colon$ Letting $1\to F\to Q\to T\to 1$ be a flasque resolution \eqref{exact_sequence:flasque-resolution}, we get the following morphism of long exact sequences \bd\label{ext-seq:flasque_resolution}\cdots\tir&H^1(A,Q)\tir\arrow[d]&H^1(A,T)\tir\arrow[d,"\theta_1^T"]&H^2(A,F)\tir\arrow[d,"\theta_2^F"]&\cdots\\\cdots\tir& H^1(K,Q)\tir&H^1(K,T)\tir&H^2(K,F)\tir&\cdots.\ed Since $Q$ is quasi-trivial, there are connected finite \'etale covers $A\to A_i$ such that $Q\cong\prod\mathrm{Res}_{A_i/A}(\bb{G}_m)$. By the fact that higher direct images vanish along finite morphisms \cite[\href{https://stacks.math.columbia.edu/tag/03QP}{Tag 03QP}]{stacks-project}, the cohomology rewrites itself as $H^1(A,Q)\cong\prod H^1(A_i,\bb{G}_m)$, and since all the rings $A_i$ are semilocal (since an \'etale morphism is quasi-finite), the cohomology vanishes thanks to the Hilbert theorem 90 \ccite{milne-etale}{Chapter III, Section 4, Proposition 4.9}. By a similar argument, the cohomology $H^1(K,Q)$ vanishes, and in view of \eqref{ext-seq:flasque_resolution}, to prove \eqref{theta1Tinjection}, it is enough to show that $\theta_2^F$ is injective. By using the fact that the \'etale cohomology commutes with filtered colimits of rings, this is a consequence of \Cref{lem:Brauer_group_injects}.\vspace{0.15cm}
				\par \eqref{theta2Tinjection}$\colon$ This time, letting $1\to T\to F'\to Q'\to 1$ be a different flasque resolution \eqref{exact_sequence:flasque-resolution2}, we get the following morphism of long exact sequences \bud\label{ext-seq:different_flasque_resolution}\cdots\tir&H^1(A,Q')\tir\arrow[d]&H^2(A,T)\tir\arrow[d, "\theta_2^T"]&H^2(A,F')\tir\arrow[d,"\theta_2^{F'}"]&\cdots\\\cdots\tir&H^1(K,Q')\tir& H^1(K,T)\tir&H^1(K,F')\tir&\cdots.\eud In a similar vein as above, since $Q'$ is quasi-trivial, it follows that $H^1(A,Q')=H^1(K,Q')=0$. Consequently, a similar argument as above shows that $\theta_2^{F'}$ is injective. Therefore, $\theta_2^T$ is injective, and we are done.
				\es
				Remarkably, for a flasque torus $F$, the $\bb{A}^1$-invariance of $H^1(-,F)$ follows as a direct consequence of \Cref{cor:torus_Gersten_injectivity}. This is encapsulated in the corollary below. This conclusion readily follows from the Noetherian case studied in \ccite{ct-sansuc_flasque_tori}{Lemma 2.4}.
				\bc\label{cor:normality}
				For a semilocal Pr\"ufer domain $R$, a semilocal, integral domain $A$ that is $R$-essentially smooth, a flasque $A$-torus and a nonempty, quasi-compact open $U\subset\bb{A}^n_A$, the composite morphism
				\bud H^1(A,F)\to H^1(\bb{A}_A^n,F)\to H^1(U,F)\text{\hspace{0.5cm}is surjective.}\eud 
				\ec
				\bs
				Thanks to \eqref{theta1F_surjection}, the second morphism is surjective. Consequently, it reduces us to proving the surjectivity of the first morphism. In fact, we shall show that  \bud H^1(A,F)\iso H^1(\bb{A}^n_A,F).\eud 
				Assuming that $A$ is normal, by a limit argument, it is enough to show the displayed isomorphism for Noetherian, normal domains, which is the content of loc.~cit. Therefore, it remains to show that $A$ is normal. This follows from \cite[\href{https://stacks.math.columbia.edu/tag/00GY}{Tag 00GY}, \href{https://stacks.math.columbia.edu/tag/030A}{Tag 030A}]{stacks-project} and \ccite{sga1}{Expos\'e I, th\'eor\`eme 9.5(i)} (cf.~\ccite{milne-etale}{Chapter I, Proposition 3.17(b)}).
				\es

\printbibliography
\end{document}